\documentclass{amsart}
\usepackage{amssymb}
\usepackage{graphicx}
\usepackage{epstopdf}
\usepackage{amscd}

\newtheorem{thm}{Theorem}[section]
\newtheorem{cor}[thm]{Corollary}
\newtheorem{lem}[thm]{Lemma}
\newtheorem{prop}[thm]{Proposition}
\newtheorem{conj}[thm]{Conjecture}
\theoremstyle{definition}
\newtheorem{defn}[thm]{Definition}
\theoremstyle{remark}
\newtheorem{rem}[thm]{Remark}

\numberwithin{equation}{section}
\newcommand{\To}{\longrightarrow}

\newcommand{\inv}{^{-1}}
\newcommand{\C}{\mathbb C}
\newcommand{\Z}{\mathbb Z}
\newcommand{\R}{\mathbb R}

\newcommand{\x}{\times}

\newcommand{\al}{\alpha}

\newcommand{\pr}{\operatorname{pr}}

 \makeatletter 
 \def \underbracket { %
 \@ifnextchar[{\@underbracket}{\@underbracket[\@bracketheight]}%
 } 
 \def\@underbracket[#1]{ %
 \@ifnextchar[{\@under@bracket[#1]}{\@under@bracket[#1][0.4em]} %
 } 
\def\@under@bracket[#1][#2]#3{
\mathop{\vtop{\m@th\ialign{##\crcr$\hfil\displaystyle{#3}\hfil$ %
\crcr\noalign{\kern3\p@\nointerlineskip}\upbracketfill
{#1}{#2} 
\crcr\noalign{\kern3\p@}}}}\limits} 
\def\upbracketfill#1#2{$\m@th\setbox\z@\hbox{$\braceld$} 
\edef\@bracketheight{\the\ht\z@}\bracketend{#1}{#2} 
\leaders\vrule\@height#1\@depth\z@\hfill
\leaders\vrule\@height#1\@depth\z@\hfill 
\bracketend{#1}{#2}$} 
\def\bracketend#1#2{\vrule height#2 width#1\relax} 
\makeatother

\newcommand{\variable}{\underbracket[.5 pt][2 pt]{\ \ }\hskip -4pt}

\begin{document}

\title
{A mirror symmetric construction of $qH^*_T(G/P)_{(q)}$}
\author{Konstanze Rietsch}%
\address{King's College London, UK}%
\email{konstanze.rietsch@kcl.ac.uk}%

\keywords{flag varieties, quantum cohomology, mirror symmetry}%

\thanks{
The work for this paper was done while the author was funded by a 
Royal Society Dorothy Hodgkin Research Fellowship. The author is
currently supported by EPSRC advanced fellowship
EP/S071395/1}

\subjclass[2000]{20G20, 15A48, 14N35, 14N15} \keywords{Flag varieties,
quantum cohomology, mirror symmetry}

\date{August 22, 2007}
\begin{abstract}
Let $G$ be a simple simply connected complex algebraic group. 
We give a Lie-theoretic construction of a conjectural mirror family 
associated to a general flag variety G/P, and show that it recovers the 
Peterson variety presentation for the $T$-equivariant 
quantum cohomology rings $qH_T^*(G/P)_{(q)}$ 
with quantum parameters inverted.
For $SL_n/B$ we relate our construction to the mirror family defined by Givental
and its 
$T$-equivariant analogue due to Joe and Kim. 
\end{abstract}
\maketitle
\section{Introduction}

According to Givental \cite{Giv:MSFlag} and Eguchi, Hori and Xiong \cite{EHX:GravQCoh}
mirror symmetry 
should have an extension to Fano manifolds $X$, where it means
in essence a `mirror side' representation of the quantum cohomology 
$D$-module, or quantum differential equations, of $X$ by complex 
oscillatory integrals. Such mirror models have previously been  
constructed for toric Fano manifolds and the flag variety
$SL_n/B$ by Givental \cite{Giv:MSFlag,Giv:MSToric}.
Moreover for $SL_n/B$  Joe and Kim proved also a $T$-equivariant
version of Givental's mirror theorem \cite{JoeKim:EquivMirrors}.  

In this paper we are interested in the case where $X$ is 
a general flag variety $G/P$. 
Explicitly the ingredients for the mirror symmetric model associated
to $X$ should be the following.
\begin{enumerate}
\item 
A $k$-parameter family $Z$ of (affine) varieties
$Z_s$ of dimension $d=\dim_\C (X)$. Here $k=\dim H^2(X,\C)$.
\item 
A family of holomorphic $d$-forms $\omega_s$ on the fibers
$Z_s$ in the family.
\item A holomorphic function $\mathcal F:Z\To \C$, which 
will play the role of the phase.
\end{enumerate}

From these data one can write down complex oscillatory integrals
\begin{equation}\label{e:S}
S_\Gamma(s)=\int_{\Gamma_s}  
e^{\mathcal F/\hbar}\ \omega_s , 
\end{equation}
where the $\Gamma$ are  certain continuous families 
of (possibly non-compact) cycles $\Gamma_s$ in $Z_s$, for example associated
to $\mathcal F$ via Morse theory of $Re(\mathcal F)$,
see \cite{Giv:MSFlag}, \cite{AGV:SingBookII}. 

In our case $X=G/P$ and has an action of a maximal torus $T$. Let 
$\mathfrak h$ be the Lie algebra of $T$. To obtain a $T$-equivariant analogue we need to 
add one more item to the data (1-3). 

 (4) A multi-valued holomorphic function $\phi:Z\x \mathfrak h\to\C$. Or more
 precisely, a holomorphic function $\tilde \phi$ on a covering $\tilde Z\x \mathfrak h$ of 
 $Z\x \mathfrak h$.

Using (4) one can write down the more general integrals
\begin{equation}\label{e:I}
\tilde S_{\Gamma}(s,h)=\int_{\Gamma_s}e^{\tilde{\mathcal F}/\hbar}\ \tilde\phi(\ ,h)\ \tilde\omega_s,
\end{equation}
where $\Gamma_s$ now lies in the covering $\tilde Z$ of $Z$ and we have denoted
the pullbacks of $\mathcal F$ and $\omega_s$ to $\tilde Z$ by $\tilde{\mathcal F}$ and 
$\tilde \omega_s$, respectively.  

Mirror symmetry for $G/P$ should involve a presentation of the set of solutions to the 
($T$-equivariant) quantum differential equations associated to $G/P$, see \cite{Giv:EquivGW,JoeKim:EquivMirrors, CoxKatz:QCohBook}, via integrals of the form \eqref{e:S},
respectively \eqref{e:I}. 

We now turn our attention to quantum cohomology. There is a remarkable, unified Lie-theoretic presentation for the ($T$-equivariant) 
quantum cohomology rings $qH_{T}^*(G/P)$ which 
was discovered by Dale Peterson \cite{Pet:QCoh}.
From his point of view the quantum cohomology rings arise 
as (posssibly non-reduced) coordinate rings,
$$
qH^*_{T}(G/P)\cong \C[\mathcal Y_P ],
$$
where $\mathcal Y_P$ is a particular affine stratum, of 
the so-called `Peterson variety'
$\mathcal Y$ in  $G^\vee/B^\vee\x \mathfrak h$. We will review 
Peterson's results in Section~\ref{s:Peterson}.

Following Givental \cite{Giv:EquivGW} on the other hand, relations for the small quantum 
cohomology ring are obtained as equations for the
characteristic variety of the quantum cohomology $D$-module.
And for this variety there is a somewhat corresponding 
construction on the mirror side, which  is to look at  
what is `swept out' by the critical points of $\mathcal F$  
along the  fibers of the family $Z$ (or the critical points of 
$\mathcal F+\ln \phi(\ ,h )$,  
in the $T$-equivariant case).

In this paper we will give a Lie theoretic construction associating 
to any $G/P$ a family $Z=Z_P$ with associated data
(1-4). The fibers $Z_s$ of the family turn out to have  natural 
compactifications to $G^\vee/P^\vee$, and the base is
the algebraic torus $H^2(G/P,\C)/2\pi iH^2(G/P,\Z)$, or is $H^2(G/P,\C)$ 
if we pull back along the exponential map.  
The main result, Theorem~\ref{t:main}, says that the critical points
of $\mathcal F+\ln \phi(\ ,h )$ along the fibers $Z_s$ and 
for varying $h$ indeed recover the Peterson variety stratum 
$\mathcal Y_{P}$ 
(or, more precisely, the open dense part
in $\mathcal Y_P$ where the quantum 
parameters are nonzero). 

This result supports the mirror conjectures, stated in Section~\ref{s:conjecture},
that the integrals
\eqref{e:S} and \eqref{e:I} defined in terms of our data 
$(Z_P, \omega,\mathcal F_P,\phi_P)$ 
should give solutions to
the quantum differential equations  
associated to $G/P$ and their $T$-equivariant analogues, respectively.

In the final section we verify these mirror conjectures in the special case of $SL_n/B$ 
by comparing our mirror construction with Givental's \cite{Giv:MSFlag},
and, in the equivariant case, with the construction of Joe and Kim, \cite{JoeKim:EquivMirrors}.
Explicitly, Givental's mirror family  
is shown to appear as an open subset inside our $Z_B$,
and $\mathcal F$ and the $\omega_s$ are related by restriction. 
The relationship with Joe and Kim's integrals is via a comparison
map which is a covering of an open inclusion, but 1-1 on any of
Joe and Kim's integration contours.   

We plan to discuss the mirror conjecture for the general $G/B$ case in a future 
paper. In that setting the quantum differential equations were
determined by Kim, and the mirror conjecture can be interpreted as 
saying that the integrals 
define Whittaker functions obeying the quantum Toda lattice
associated to the Langlands dual root system.   In this direction,
but still confined to type $A$, there
has already been some interesting independent work of Gerasimov, Kharchev,
Lebedev and Oblezin \cite{GKLO:GaussGivental}, who reproved Givental's mirror theorem
using representation theory.

This work was motivated on the one hand by a desire to 
put the papers \cite{Giv:MSFlag,JoeKim:EquivMirrors,BCKS:MSPFlag}, 
concerning mirror constructions for classical flag varieties, 
into a Lie theoretic context.
And on the other hand it is an attempt to better understand 
Dale Peterson's powerful point of view about quantum
cohomology, \cite{Pet:QCoh}. Morally speaking,
it says that Peterson's presentation of $qH_T^*(G/P)$
via the variety $\mathcal Y_P$ might
be considered 
a mirror symmetry phenomenon for $G/P$. 

\vskip .2cm

\noindent{\it Acknowledgments.} I am indebted to Dale Peterson for his 
inspiring lectures on quantum cohomology. These results were finally 
written up during a six month stay in Waterloo, Canada. I would like to thank
the Perimeter Institute and the University of Waterloo for their hospitality.

\section{Background and notation}

Let $G$ be a simple simply connected algebraic group over $\C$
of rank $n$.  
We fix opposite Borel subgroups $B=B_-$ and $B_+$ with unipotent radicals 
$U_-$ and $U_+$, respectively. Let $T$ be 
the maximal torus $T=B_+\cap B_-$, and 
$W=N_G(T)/T$ the Weyl group.

Let $\mathfrak g$ be the Lie algebra of $G$ and  
$\mathfrak b_-,\mathfrak b_+,\mathfrak u_-,\mathfrak u_+, 
\mathfrak h$ the Lie algebras of $B_-,B_+,U_-$,$U_+$ and $T$, respectively.
The adjoint action of $G$ is denoted by a dot for simplicity. 
So $g\cdot X:=\operatorname{Ad}(g)X$, for $g\in G$ and $X\in\mathfrak g$.
Similarly for the coadjoint action, so when $X\in\mathfrak g^*$.

Let $X^*(T)$ be the character group of $T$ and $Q\subset X^*(T)$ the root
lattice. We will sometimes view these as lying in $\mathfrak h^*$. Let  $\Delta_+$ be the set of positive roots
corresponding to $B_+$, so that the Lie algebra of $B_+$
written as sum of weight spaces with respect to the adjoint action of $T$ is
$$
\mathfrak b_+=\bigoplus_{\alpha\in\Delta_+}\mathfrak g_{\alpha}.
$$ 
We set $I=\{1,\dotsc,n \}$, where $n$ is the rank of $G$, and 
use $I$ to enumerate the simple roots $\{\alpha_i \ | \ i\in I\}$ in $\Delta_+$. 
Corresponding to the simple roots (and their negatives), we have the Chevalley generators 
$e_i$ and $f_i$ in $\mathfrak g_{\alpha_i}$ and $\mathfrak g_{-\alpha_i}$, respectively. 
These define the one parameter subgroups 
 \begin{equation*}
x_i(t):=\exp(t e_i), \qquad y_i(t):=\exp(t f_i),
 \end{equation*}
where $t\in\C$.
Let 
\begin{equation}\label{e:Weyl}
\dot s_i=x_i(1)y_i(-1)x_i(1).
\end{equation}
Then $\dot s_i$ represents a simple
reflection in $W$ which we denote by $s_i$.
For general $w\in W$, a representative
$\dot w\in G$ is defined by $\dot w=\dot s_{i_1}\dot
s_{i_2}\cdots \dot s_{i_m}$, where $ s_{i_1} s_{i_2}\cdots
s_{i_m}$ is a (any) reduced expression for $w$. The length $m$ of
a reduced expression for $w$ is denoted by $\ell(w)$.

Let $P\supseteq B$ be a (fixed) parabolic subgroup of $G$. 
Define
$I_P=\{i\in I\ |\ \dot s_i\in P \}$ and let $I^P$ be its complement in
$I$. We will usually denote the elements of $I^P$ by
\begin{equation*}
I^P=\{n_1,\dotsc, n_k\},
\end{equation*}
for $1\le n_1<n_2<\cdots<n_k$. We denote by $W_P$ the parabolic subgroup of $W$ associated to 
$P$, and by $W^P$ the set 
of minimal length coset representatives in $W/W_P$. So
 \begin{align*}
 W_P&:=\left<s_i\ |\ i\in I_P\right >, \\
 W^P&:=\{w\in W\ |\ \ell(w s_i)>\ell(w) \text{ for all $i\in I_P$}.
 \}
 \end{align*}
Let $w_P$ be the longest element in the parabolic subgroup $W_P$. 
For example $w_B=1$ and $w_G$ is the longest element in $W$, also 
denoted $w_0$.

Let $G^\vee$ be the Langlands dual group to $G$. Note that $G^\vee$ is adjoint
since $G$ was simply connected. We will use all the same notation for $G^\vee$
as for $G$, but with an added superscript where required. 
For example the Chevalley generators 
of $\mathfrak g^\vee$ are denoted by $e_i^\vee$ and $f_i^\vee$, where $i\in I$. 
The Weyl group for $G^\vee$ is again $W$. For simplicity we will write $\dot w$ 
again
for the representative of $w$ in $G^\vee$ obtained as above. 
Identify $\mathfrak h^\vee$ with $\mathfrak h^*$, the dual of the Lie algebra of $T$.
In particular we may view the weight and root lattices of $G$ as lying inside $\mathfrak h^\vee$.
The dual pairing between $\mathfrak h$ and $\mathfrak h^\vee$ is denoted by $<\ ,\ >$. 

We will also consider the universal covering group $\widetilde G^\vee$ 
of $G^\vee$. Let $\pi~:\widetilde G^\vee\To G^\vee$ be the covering map. 
The group $\widetilde G^\vee$ has maximal torus $\tilde T^\vee=\pi\inv(T^\vee)$ 
and Borel 
subgroups $\tilde B_-^\vee=\pi\inv(B^\vee_-)$ and
$\tilde B_+^\vee=\pi\inv(B^\vee_+)$. The unipotent 
radicals of $\tilde B^\vee_-$ and $\tilde B^\vee_+$
can be identified with $U^\vee_-$ and $U^\vee_+$,
respectively, via $\pi$. Also the Weyl group representatives 
in $\widetilde G^\vee$ defined via \eqref{e:Weyl} are 
identified via $\pi$ with those in $G^\vee$, and we 
will suppress the difference in our notation.

For any dominant coweight $\lambda^\vee$ we have an irreducible 
representation $V(\lambda^\vee)$ of $\widetilde G^\vee$. In each $V(\lambda^\vee)$ 
let us fix a lowest weight vector $v^-_{\lambda^\vee} $. 
Then for 
any $v\in V(\lambda^{\vee})$ and extremal weight vector $\dot w\cdot v^-_{\lambda^\vee}$ 
we have the coefficient $\left<v,\dot w\cdot v^-_{\lambda^\vee}\right>\in \C$ defined by
\begin{equation*}
v=\left<v,\dot w\cdot v^-_{\lambda^\vee}\right>\dot w\cdot  v^-_{\lambda^\vee}+
 \text{other weight space summands.}
\end{equation*}
Let  $v^+_{\lambda^\vee}:=\dot w_0\cdot v^-_{\lambda^\vee}$. The most important
choices for $\lambda^\vee$ are the fundamental coweights $\omega_i^\vee$, where $i\in I$, 
and $\rho^\vee:=\sum_{i\in I}\omega^\vee_i$. If $\lambda^\vee\in Q^\vee$ then also
$G^\vee$ acts on $V(\lambda^\vee)$.

In the Langlands dual context we will consider the flag variety $G^\vee/B^\vee_-$.  Then 
for two elements $v,w\in W$ with $v\le w$ we have the intersection of opposed
Bruhat cells
\begin{equation*}
\mathcal R^\vee_{v,w}:=(B^\vee_+\dot v B^\vee_-
\cap B^\vee_-\dot w B^\vee_-)/B^\vee_-
\end{equation*}
in $G^\vee/B^\vee_-$. It is known that $\mathcal R^\vee_{v,w}$ is smooth and
irreducible of dimension $\ell(w)-\ell(v)$, \cite{KaLus:Hecke2,Lus:IntroTotPos}.  

\section{Equivariant quantum cohomology of $G/P$ and Peterson's presentations.}
\subsection{}
The (small) quantum cohomology ring of $G/P$ is a deformation of the usual 
cohomology ring with $k=\dim H^2(G/P)$ parameters,
\begin{equation*}
qH^*(G/P)\cong H^*(G/P)\otimes \C[q_1,\dotsc, q_k],
\end{equation*}
where the deformed cup product has structure constants given by 
genus 0, 3-point Gromov-Witten invariants.
We refer the reader to \cite{CoxKatz:QCohBook,FuPa:QCoh} for definitions and
background. 
Note that we will always take
coefficients to be in $\C$.
For an equivariant version of quantum cohomology  
see the papers \cite{Beh:LocalizationGW, GiKi:FlTod, Kim:EquivQCoh}. 
The $T$-equivariant quantum cohomology $qH^*_T(G/P)$ is a 
module over $\C[q_1,\dotsc, q_k]$ and $\C[\mathfrak h]$, and 
is a simultaneous deformation of the quantum cohomology and the equivariant cohomology rings. 

In the literature there are many special cases of flag varieties where presentations of 
quantum  cohomology rings have been explicitly determined. See for example 
\cite{Kim:QCohPFl,AstSa:QCohPFl,Cio:QCohPFl,Kim:EquivQCoh} in type $A$, \cite{Kim:QCohG/B} for general $G/B$, 
and \cite{Tamvakis:QCohGrass} for Grassmannians in other types.

 The structure of (non-equivariant) quantum cohomology for general $G/P$ 
is described in \cite{FuWo:SchubProds,Wood:Compare,Pet:QCoh}.  
For  $qH^*_T(G/P)$ Mihalcea has given a quantum Chevalley formula
\cite{Mih:EquivQCoh}, and thereby completely determined the  
ring structure. 

The only general construction of presentations for quantum cohomology 
rings of flag varieties $G/P$ is due to Dale Peterson \cite{Pet:QCoh}, unpublished
so far. It involves the 
remarkable `Peterson variety' $\mathcal Y$ which we now introduce.

\subsection{} \label{s:Peterson}
Following \cite{Pet:QCoh} we define a closed $2n$-dimensional
subvariety  
$\mathcal Y$ of $G^\vee/B^\vee\x \mathfrak h\,$.
Let us canonically identify $\mathfrak h$ with 
the zero weight space $(\mathfrak g^\vee)^*$  
via $\mathfrak h\cong (\mathfrak h^\vee)^*$.  
Define 
$$
F:=\sum_{i\in I}(e_i^\vee)^*\ \ \in \ \ (\mathfrak g^\vee)^*,
$$
where $(e_i^\vee)^*$ denotes the 
linear functional which is one on $e_i^\vee$ and zero along all other weight 
spaces. We write $g\cdot \eta$ for the coadjoint action of $g\in G^\vee$ on 
$\eta\in\mathfrak (g^\vee)^*$. The (equivariant) Peterson variety is the subvariety 
of $G^\vee/B_-^\vee\x\mathfrak h$ defined by
\begin{equation*}
\mathcal Y:=
\left\{\ (gB^\vee_-,h)\in (G^\vee/B^\vee_-)\x\mathfrak h 
\ |\  g\inv\cdot (F-h) \text{ vanishes on } [\mathfrak u^\vee_-,\mathfrak u^\vee_-] \right \}.
\end{equation*}
Its fiber over $0\in\mathfrak h$ is
$$
Y:=\left \{\ gB^\vee_-\ |\ g\inv\cdot F\text{ vanishes on } 
[\mathfrak u^\vee_-,\mathfrak u^\vee_-] \right \},
$$
and may also be called the Peterson variety, see  \cite{Kos:QCoh}. 

To a parabolic $P\supseteq B$ associate strata
$\mathcal Y_P$ and $Y_P$, in $\mathcal Y$ and $Y$, respectively, 
which arise from (possibly non-reduced) intersections
with Bruhat cells for $B^\vee_+$,  
\begin{align*}
\mathcal Y_P &:=\mathcal Y\x_{(G^\vee/B^\vee_-)\x\mathfrak h} 
(B_+^\vee \dot w_P B^\vee_-/B^\vee_-\x \mathfrak h),\\
Y_P &:=Y\x_{G^\vee/B^\vee_-} 
(B_+^\vee \dot w_P B^\vee_-/B^\vee_-).
\end{align*}
Moreover we consider the following open subvarieties obtained by intersection with 
the big Bruhat cell for $B^\vee_-$,
\begin{align*}
\mathcal Y^*:=&\mathcal Y\x_{(G^\vee/B^\vee_-)\x\mathfrak h}
(B^\vee_- \dot w_0 B^\vee_-/B^\vee_-\x\mathfrak h), \\
Y^*:=&Y\x_{G^\vee/B^\vee_-} (B^\vee_-\dot w_0 B^\vee_-/B^\vee_-),
\end{align*}
and their strata
\begin{align*}
\mathcal Y_P^*&:=\mathcal Y\x_{(G^\vee/B^\vee_-)\x\mathfrak h}
(\mathcal R^\vee_{w_P,w_0}\x \mathfrak h),\\
Y_P^*&:=\mathcal Y\x_{G^\vee/B^\vee}
\mathcal R^\vee_{w_P,w_0}.
\end{align*} 

\subsection{}
We now state some results of Peterson's \cite{Pet:QCoh} which are essential
to, and significantly inspired, this work.

\subsubsection{} First of all,
$\mathcal Y_{P}$ and $Y_P$ are (possibly non-reduced) affine varieties
of pure dimension $|I^P|+ n$ and $|I^P|$, respectively, and one has the decomposition
\begin{align*}
\mathcal Y(\C)&=\bigsqcup_P\mathcal Y_P(\C)
\end{align*}
for the $\C$-valued points. Here $P$ runs over the set of all
parabolic subgroups of $G$ containing $B$. 
See \cite{Kos:QCoh, Kos:QCoh2} for a treatment of 
the non-equivariant case, in particular Kostant proved 
that $Y_B$ is irreducible.   

\subsubsection{}\label{s:stabilizer} There is an isomorphism, 
\begin{equation*}
\mathcal Y^*\cong\left\{ (b,h)\in B^\vee_-\x\mathfrak h\ |\ 
b\cdot (F-h)=F-h
\right\}
\end{equation*}
given by the map $(b,h)\mapsto (b\dot w_0 B^\vee_-,h)$ from right to left. 
Setting $h=0$ we have that 
$Y^*$ is isomorphic to the stabilizer in $B^\vee_-$ of $F$.

\subsubsection{}
There is an explicit isomorphism
\begin{equation}\label{e:P-iso1}
qH^*(G/P)\overset\sim\To \C[Y_P],
\end{equation}
from the quantum cohomology ring of $G/P$
to the coordinate ring of $Y_P$.   
See \cite{Rie:QCohPFl, Rie:ErratumJAMS} for a description and 
proof of Peterson's isomorphism in the case
 of $qH^*(SL_n/P)$.

Replacing $Y$ by $\mathcal Y$ and $qH^*(G/P)$ by $qH^*_T(G/P)$
in \eqref{e:P-iso1} gives an equivariant version of this result, 
\begin{equation}\label{e:P-isoT}
qH_T^*(G/P)\overset\sim\To \C[\mathcal Y_P],
\end{equation}
which was formulated
in \cite{Pet:QCoh}, and follows from \cite{Pet:QCoh} and \cite{Mih:EquivQCoh}. 

\subsubsection{}The varieties $\mathcal Y_P^*, Y_P^*$ are 
open dense in 
$\mathcal Y_P$ and $Y_P$, respectively,  
and the map \eqref{e:P-isoT} induces an isomorphism
\begin{equation}\label{e:P-iso2T}
qH_T^*(G/P)[q_1\inv,\dotsc,q_{k}\inv]\overset\sim\To 
\C[\mathcal Y_P^*].
\end{equation}

\subsubsection{}\label{s:G/B}
In the case of $G/B$ the isomorphism 
\begin{equation}\label{e:G/Biso}
qH_T^*(G/B)\overset\sim\To \C[\mathcal Y_B],
\end{equation}
is related to Kim's presentation \cite{Kim:QCohG/B}  of $qH^*(G/B)$ as follows.
Kim described the relations of the $T$-equivariant small quantum cohomology
ring of $G/B$ in terms of integrals of motion of the Toda lattice
associated to the Langlands dual group. The phase
space $T^*(T^\vee)\hat=T^\vee\x\mathfrak h$ 
of the Toda lattice for $ G^\vee$ may be embedded into $(\mathfrak g^{\vee})^*$ by
$$
(t,h')\mapsto F- h' -\sum_{i\in I} \alpha_i^\vee(t) {(f_i^\vee)}^* ,
$$
where ${(f_i^\vee)}^*$ is defined analogously to ${(e_i^\vee)}^*$ and 
$\mathfrak h$ is identified with $(\mathfrak h^\vee)^*$ viewed as a subspace of 
$(\mathfrak g^\vee)^*$. This is Kostant's construction  \cite{Kos:Toda}. The 
image of the embedding is the translate by $F$
of a $B_-^\vee$-coadjoint orbit in $(\mathfrak b^\vee_-)^*\subset(\mathfrak g^\vee)^*$,  
and the integrals of motion of the Toda lattice are 
given by restrictions of 
$G^\vee$-invariant polynomials on $(\mathfrak g^\vee)^*$.
By Chevalley's restriction theorem  $\C[(\mathfrak g^\vee)^*]^{G^\vee}=\C[\mathfrak h]^W$
and the latter is a polynomial ring with $n$ homogeneous generators $\Sigma_1,\dotsc, \Sigma_n$. 
Now let 
$$
\mathcal  A:= F + \mathfrak h\oplus \sum\C {(f_i^\vee)}^* \subset (\mathfrak g^\vee)^*,
$$
and consider the map 
$$
\Sigma:\mathcal A \to \mathfrak h/W
$$
 obtained from $\C[\mathfrak h]^W\hookrightarrow \C[(\mathfrak g^\vee)^*]$.
Then Kim's presentation of $qH^*_T(G/B)$ can be restated as an isomorphism 
$$
qH^*_T(G/B)\overset\sim\to \C[\mathcal A\x_{\mathfrak h/W}\mathfrak h]. 
$$
Finally, we have a map
\begin{eqnarray*}
\mu:\mathcal Y_B&\to & \mathcal A\x_{\mathfrak h/W}\mathfrak h\\
(uB^\vee_-,h) &\mapsto & (u\inv\cdot(F-h),h),
\end{eqnarray*}
where $u\in U^\vee_+$. This is an isomorphism by another result of Kostant's,
see \cite{Kos:PolReps}. 
Peterson's map \eqref{e:G/Biso} is given by the composition of Kim's presentation with $\mu^*$. 
\subsubsection{} For $w\in W^P$ let 
$$
\sigma^w_{G/P}\in qH^*_T(G/P)
$$ 
denote the corresponding (quantum equivariant) 
Schubert class. The Schubert classes $\sigma^w_{G/B}$ for $G/B$ may be viewed as 
rational functions on $\mathcal Y$. Peterson's theory implies that  
if $w\in W^P$, then the restriction of $\sigma^w_{G/B}$
to $\mathcal Y_P$ is a regular function and 
under \eqref{e:P-isoT} represents the Schubert class $\sigma^w_{G/P}$. 
In particular it follows that all of the isomorphisms \eqref{e:P-isoT}
for varying $P$ are explicitly determined by \eqref{e:G/Biso}. 
In the special case of $qH^*(SL_n/P)$ 
an ad hoc proof of this relationship between the Schubert classes
is given in \cite{Rie:QCohPFl}.

\subsubsection{}\label{s:qparams}
Let $j\in\{1,\dotsc, k\}$. The map \eqref{e:P-isoT} 
identifies $q_j$ with the regular function on $\mathcal Y_P$
given by
\begin{eqnarray*}
(u\dot w_P B^\vee_-,h)&\mapsto & -( F-h)(u\dot w_P\cdot f_{n_j}^\vee).
\end{eqnarray*}

\section{A mirror construction for $H^*_T(G/P)_{(q)}$}\label{s:mirrors}
In this section we will introduce the ingredients (1), (3) and (4) of mirror symmetry 
described in the  
introduction for a general flag variety $G/P$. Then we will state the main theorem,
which gives a mirror symmetric construction of the strata 
$\mathcal Y_P^*$ in the equivariant Peterson variety.     
\subsection{}\label{s:mirrors}
Let
\begin{equation}\label{e:Z}
Z=Z_P:=
\{(t,b)\in {(T^\vee)}^{W_P}\x B^\vee_-\ 
|\ b\in U^\vee_+ t \dot w_P\dot w_0\inv U^\vee_+\}.
\end{equation}
We view $Z_P$ as a family of varieties via the map 
$\pr_1: Z_P\to {(T^\vee)}^{W_P}$
projecting onto the first factor. For $t\in (T^\vee)^{W_P}$ let us  write 
\begin{equation}\label{e:fiber}
Z_P^t:= B^\vee_-\cap U^\vee_+ t \dot w_P\dot w_0\inv U^\vee_+,
\end{equation}
which we may identify with the fiber $\pr_1\inv(t)$ in $Z_P$. 
We record the following basic properties of  the family $Z_P$. 
\begin{enumerate}
\item
Projection onto the second factor in $Z_P$ restricts to an isomorphism
 $$\pr_2:Z_P\overset\sim\To B^\vee_-\cap U^\vee_+(T^\vee)^{W_P}
 \dot w_P\dot w_0\inv U^\vee_+.$$ 
\item
Fix $t\in (T^\vee)^{W_P}$. Then the fiber $Z^t_P$ is smooth of 
dimension $n_P=\dim G/P$, and may be identified with $\mathcal R^\vee_{w_P,w_0}$
by
\begin{eqnarray*}
Z^t_P&\To &\mathcal R^\vee_{w_P,w_0}\\
b &\mapsto & b\dot w_0 B^\vee_-/B^\vee_-.
\end{eqnarray*}
\item
Using the isomorphism from 
(2) to identify all the fibers we obtain a trivialization
\begin{equation}\label{e:triv}
Z_P\overset\sim\To (T^\vee)^{W_P}\x\mathcal R^\vee_{w_P,w_0}.
\end{equation}
The map $\psi_P: Z_P\to \mathcal R^\vee_{w_P,w_0}$ obtained by 
composing with the projection onto the second factor in the trivialization
will be important later on. 
\end{enumerate}
The properties (1-3) are straightforward to verify. 
We note that the fibers can be naturally compactified to give the 
Langlands dual flag variety $G^\vee/P^\vee$. 

\subsection{}\label{s:F} 
Let 
\begin{eqnarray*}
f^\vee&=&\sum_{i\in I} f^\vee_i, \\
f^\vee_{(P)}&=& \sum_{i\in I}\frac{1}{<\rho^\vee,w_P\cdot\alpha_i>} f^\vee_i.
\end{eqnarray*}
In particular $ f^\vee_{(B)}=f^\vee$. We define a function $\mathcal F=
\mathcal F_P:Z_P\to \C$ in terms of the representation $V(\rho^\vee)$ of $\widetilde G^\vee$
as follows.
\begin{equation}\label{e:F}
\mathcal F_P(t,b)=
\frac{\left<f^\vee_{(P)}\tilde b\cdot v^+_{\rho^\vee},\dot w_P\cdot 
v^-_{\rho^\vee}\right>+
\left<\tilde b f^\vee\cdot v^+_{\rho^\vee},
\dot w_P\cdot v^-_{\rho^\vee}\right>}{\left<\tilde b\cdot v^+_{\rho^\vee},
\dot w_P\cdot v^-_{\rho^\vee}\right>}
\end{equation}
where $\tilde b\in \widetilde G^\vee$ with 
$\pi(\tilde b)=b$. Note that the denominator 
insures that $\mathcal F_P$ is well defined, that is, independent of the choice 
of $v^-_{\rho^\vee}$ or lift $\tilde b$. We also denote by  $\mathcal F_P$ 
the restriction to any  $Z_P^t$. 

\subsection{}\label{s:phi}
Consider the fundamental representations $V(\omega_i^\vee)$.
In terms similar to \eqref{e:F}  the multi-valued function
$\phi=\phi_P:Z_P\x\mathfrak h\to\C$ we will define can 
be thought of as taking the form
\begin{equation}\label{e:phi}
\phi(t,b;h)=\prod_{i\in I}
\left<\tilde b\cdot v^+_{\omega^\vee_i},
v^+_{\omega_i^\vee} \right>^{\alpha_i(h)}
\end{equation}
This is only a well-defined function if $h$, after identifying $\mathfrak h$ with
$(\mathfrak h^\vee)^*$, is in the root lattice for $G^\vee$. If $h$ is also dominant we
can simply look at the representation $V(h)$, and $\phi$ becomes the
highest weight coefficient
$$
\phi(t,b;h)=\left<b\cdot v^+_h,v^+_h\right>.
$$
To give the definition more generally we consider the covering space
$$
\tilde Z_P:=Z_P\x_{T^\vee}\mathfrak h^\vee=\{(t,b,h^\vee_R)\in Z_P\x\mathfrak h^\vee\ |\ 
b\exp(-h^\vee_R)\in U^\vee_-\}
$$
of $Z_P$. If $t\in (T^\vee)^{W_P}$, let us also write
\begin{equation*}
\tilde Z^t_P:=\{(b,h^\vee_R)\in Z^t_P\x\mathfrak h^\vee\ |\ 
b\exp(-h^\vee_R)\in U^\vee_-\},
\end{equation*}
in correspondence with \eqref{e:fiber}. Then we have two families
of varieties
related by a covering map $c_P$,
\begin{align*}
\tilde Z_P\ &\overset{c_P}\To  Z_P \\
\pr_1\downarrow\ & \qquad \downarrow{\pr_1}\\
(T^{\vee})^{W_P}& = \ (T^{\vee})^{W_P},
\end{align*}
where $\tilde Z^t_P$ naturally identifies with a fiber on the left hand side, and
such that each of these fibers is also a covering of the corresponding 
fiber of $Z_P$,
\begin{equation*}
\tilde Z^t_P \to Z^t_P\ : \ (b,h^\vee_R)\mapsto  b.
\end{equation*}

We define a holomorphic function $\tilde \phi$ on $\tilde Z_P\x\mathfrak h$
by
\begin{equation}\label{e:tildephi}
\begin{array}{lcl}
\tilde \phi:\tilde Z_P\x\mathfrak h &\to &\C\\
(t,b,h^\vee_R;h)&\mapsto&e^{<h,h^\vee_R>},
\end{array}
\end{equation}
where $<\ ,\  >$ is the dual pairing between $\mathfrak h$ and $\mathfrak h^\vee$.
It is clear that this agrees with the matrix coefficient $\left<b\cdot v^+_h,v^+_h\right>$
if $h$ is a dominant weight in the root lattice of $G^\vee$. We denote the restriction of $\tilde \phi$ to any 
$\tilde Z_P^t$ again by $\tilde \phi$.

We now take \eqref{e:tildephi} to be our
definition of the multi-valued function $\phi$. 
While $\phi$ is multi-valued on $Z_P$, note that it follows immediately from 
the definition that 
the logarithmic derivative of $\tilde \phi$ along any $\tilde Z_P$ direction is 
independent of the chosen branch (i.e. depends only on $(t,b)$ and not on $h^\vee_R$).
In particular for fixed $h\in\mathfrak h$ it makes sense to talk about
critical points of $\ln(\phi (\ ;h))$ in a fiber $Z^t_P$ of the original mirror family $Z_P$, 
as we will do below. 

\subsection{}\label{s:maintheorem}
We can now formulate our main result connecting the  mirror data constructed above
with the quantum cohomology rings of the homogeneous spaces $G/P$. 
Let 
$$
Z_{P,T}^{crit}:=\{(t,b;h)\in Z_P\x\mathfrak h\ |\ 
\text{$b$ is a critical point for }\left( 
\mathcal F_P+\ln \phi(\ ;h)\right)|_{Z^t_P}\}.
$$

Note that the quantum parameters 
$q_1,\dotsc, q_k$  in $qH^*(G/P)$ can be naturally thought of as functions $e^{t_j}$
on $H^2(G/P,\C)/2\pi i H^2(G/P, \Z)$, where the $t_j$ run through a certain basis in 
$H_2(G/P)$ (dual to the Schubert basis of $H^2(G/P)$). If we identify
$$
H^2(G/P,\C)= (\mathfrak h^\vee)^{W_P}
$$ 
by the Borel-Weil homomorphism then the $t_j$ are represented by the roots $\alpha_{n_j}^\vee$ 
(of $G^\vee$) associated to $I^P$. Therefore the
$q_j$ are identified with the corresponding functions on 
${(T^\vee)}^{W_P}$, which we  again denote
by $\alpha_{n_j}^\vee$. This is precisely how the 
quantum parameters will appear below.

\begin{thm}\label{t:main}
The map
$\psi_P:Z_P\to \mathcal R^\vee_{w_P,w_0}$ from (3) in Section~\ref{s:mirrors}  
induces  an isomorphism 
\begin{equation}\label{e:main}
\psi_P\x id_{\mathfrak h}:~ Z_{P,T}^{crit}\overset \sim\To \mathcal Y_P^*,
\end{equation}
such that the following diagram commutes
\begin{equation}\label{e:diagram}
\begin{CD} Z_{P,T}^{crit} &@>{\sim} >> &\mathcal Y_P^*\\
@V{\pr_1 }VV & &@VV{(q_i)_{i=1}^k }V\\
(T^\vee)^{W_P}&@>{\sim}>> &(\C^*)^k.
\end{CD}
\end{equation}
Here the isomorphism $ (T^\vee)^{W_P}\overset\sim\to (\C^*)^k$
is given by $(\alpha_{n_j}^\vee)_{j=1}^k$. Moreover we have
\begin{equation}\label{e:stabilizer}
Z^{crit}_{P,T} =\{(t,b;h)\in Z_P\ |\ b\cdot (F-h)=F-h \}.
\end{equation}
\end{thm}

\begin{cor}\label{s:maincor} Combining \eqref{e:main} with 
the isomorphism \eqref{e:P-iso2T} one obtains
\begin{equation*}
qH^*_T(G/P)[q_1\inv,\dotsc,q_k\inv]\overset\sim\To \C[Z_{P,T}^{crit}].
\end{equation*}
\end{cor}

\section{Proof of Theorem~\ref{t:main}}
We prove first some preparatory lemmas.

\begin{lem}\label{l:constants} Let $i\in I$ and $i^*$ be such that $w_0\cdot \alpha_i=-\alpha_{i^*}$. 
\begin{enumerate}
\item Then
\begin{equation*}
\dot w_0\cdot f_i^\vee=\dot w_0\inv\cdot f_{i}^\vee=-e_{i^*}^\vee.
\end{equation*} 
\item For any parabolic $P$,
\begin{equation*}
\dot w_P\dot w_P\in (T^\vee)^{W_P}.
\end{equation*}
\item
If $i^*$ lies in $I_P$, and $\bar i\in I_P$ is defined by  $w_Pw_0\inv\cdot \alpha_i=\alpha_{\bar i}$, then
\begin{equation*}
\dot w_P\dot w_0\inv\cdot f^\vee_i=f^\vee_{\bar i}.
\end{equation*} 
\end{enumerate}
\end{lem}

\begin{proof}
Note that we have $\dot s_i\inv\cdot f_i^\vee=
-e_i^\vee$, which can be checked by a direct calculation. Similarly 
$\dot s_i\inv\cdot e_i^\vee=-f_i^\vee$.   
Now consider the fundamental representation $V(\omega_i^\vee)$ of $\widetilde G^\vee$.
We have
\begin{multline*}
\left<(\dot w_0\inv\cdot f_i^\vee)\cdot v^-_{\omega_i^\vee},
 e_{i^*}^\vee\cdot v^-_{\omega_i^\vee}\right>
 =
 \left<\dot w_0\inv f_i^\vee\cdot v^+_{\omega_i^\vee},
 e_{i^*}^\vee\cdot v^-_{\omega_i^\vee}\right>
 \\
 =
 -\left<\dot w_0\inv \dot s_i e_i^\vee {\dot s_i}\inv\cdot v^+_{\omega_i^\vee},
 e_{i^*}^\vee\cdot v^-_{\omega_i^\vee}\right>
 =-\left<\dot s_{i^*}\dot w_0\inv e_i^\vee f_i^\vee\cdot v^+_{\omega_i^\vee},
 e_{i^*}^\vee\cdot v^-_{\omega_i^\vee}\right>
 \\
 =-\left<\dot s_{i^*}\dot w_0\inv \cdot v^+_{\omega_i^\vee},
 e_{i^*}^\vee\cdot v^-_{\omega_i^\vee}\right>
=-\left<\dot s_{i^*}v^-_{\omega_i^\vee},
 e_{i^*}^\vee\cdot v^-_{\omega_i^\vee}\right>=-1.
 \end{multline*}
This implies the second equality in (1). Analogously we can show the identity 
$$
\left<(\dot w_0\cdot f^\vee_i)\cdot v^-_{\omega_i^\vee},
e_{i^*}^\vee\cdot v^-_{\omega_i^\vee} \right>=-1,
$$
and this implies also the first equality.

For (2) let $\epsilon=\dot w_P\dot w_P$. Then $\epsilon\in T^\vee$, 
and we need to show that  $\alpha^\vee_i(\epsilon)=1$ whenever $i\in I_P$. This holds since
by (1) we have
\begin{equation*}
\dot w_P\dot w_P\cdot e^\vee_{i}=e_i^\vee, \qquad\text{for $i\in I_P$.}
\end{equation*}

Applying (1) twice as follows,
\begin{equation*}
\dot w_P\inv\cdot f_{\bar i}^\vee= - e_{i^*}^\vee= \dot w_0\inv\cdot f^\vee_{i}, \qquad \text{for $i^*
\in I_P$,}
\end{equation*}   
implies (3).
\end{proof}
\begin{lem}\label{l:formula}
Let
$$b\in B^\vee_-\ \cap\ U^\vee_+(T^\vee)^{W_P}\dot w_P\dot w_0\inv U^\vee_+
$$
with factorization $b=u_1  t\dot w_P\dot w_0\inv u_2\inv$ for $u_1,u_2\in U^\vee_+$
and $t\in (T^\vee)^{W_P}$.  
Then
\begin{equation*}
\mathcal F_P(b)=
F(u_2\cdot \rho)-F(u_1\cdot \rho).
\end{equation*}
\end{lem}
\begin{proof} Let $\tilde
t \in (\widetilde T^\vee)^{W_P}$ with $\pi(\tilde t)=t$ and  $\tilde b=u_1\tilde t
\dot w_P\dot w_0\inv u_2\inv\in \widetilde G^\vee$ covering $b$. Note that 
$\left<\tilde b\cdot v^+_{\rho^\vee},\dot w_P\cdot v^-_{\rho^\vee}\right>=
\left<u_1\tilde t\dot w_P\cdot v^-_{\rho^\vee},\dot w_P\cdot v^-_{\rho^\vee}\right >=
\rho^\vee(\tilde t\, )\inv$. Then we have
\begin{multline}\label{e:sum}
\mathcal F_P(b)=
\frac{1}{\rho^\vee(\tilde t\,)\inv}
\left (\left<f_{(P)}^\vee \tilde b\cdot v^+_{\rho^\vee},
\dot w_P\cdot v^-_{\rho^\vee}\right> +
\left<\tilde b f^\vee\cdot v^+_{\rho^\vee},\dot w_P\cdot 
v^-_{\rho^\vee}\right>\right )\\
=\rho^\vee(\tilde t\,)
\left (\left<f_{(P)}^\vee u_1\tilde t\dot w_P \dot w_0\inv\cdot v^+_{\rho^\vee},
\dot w_P\cdot v^-_{\rho^\vee}\right> +
\left<u_1\tilde t \dot w_P\dot w_0\inv u_2\inv f^\vee\cdot v^+_{\rho^\vee},\dot w_P\cdot 
v^-_{\rho^\vee}\right>\right )\\=
\left<f_{(P)}^\vee u_1\dot w_P\cdot v^-_{\rho^\vee},
\dot w_P\cdot v^-_{\rho^\vee}\right> +
\left<u_1\dot w_P\dot w_0\inv u_2\inv f^\vee\cdot v^+_{\rho^\vee},\dot w_P\cdot 
v^-_{\rho^\vee}\right>.
\end{multline}
Here the $\rho^\vee(\tilde t\,)$ was cancelled against the $\tilde t$ factors in both summands. 
Note that 
\begin{equation}\label{e:weightspaces}
s_iw_P\cdot (-\rho^\vee)-w_P\cdot(-\rho^\vee)
\in \begin{cases}\Z_{<0}\,\alpha_i^\vee &
\quad \text{ for $i\in I_P$,}\\
\Z_{>0}\,\alpha_i^\vee&\quad \text{ for $i\in I^P$.}
\end{cases}
\end{equation}
Therefore if $i\in I_P$ then
$ \dot w_P\cdot v^-_{\rho^\vee}$ is annihilated by $e_i^\vee$, and if $i\in I^P$ then 
it is annihilated by $f^\vee_i$. Now the left hand summand of \eqref{e:sum}
simplifies to
\begin{multline*}
\left<f_{(P)}^\vee u_1\dot w_P\cdot v^-_{\rho^\vee},
\dot w_P\cdot v^-_{\rho^\vee}\right>=
\sum_{i\in I}\frac{1}{<w_P\cdot\rho^\vee,\alpha_i>}
\left<f_i^\vee u_1\dot w_P\cdot v^-_{\rho^\vee},\dot w_P\cdot v^-_{\rho^\vee} \right>\\
=\sum_{i\in I}\frac{1}{<w_P\cdot\rho^\vee,\alpha_i>}
(e_i^\vee)^*(u_1\cdot \rho)\left< f_i^\vee e_i^\vee\dot w_P\cdot v^-_{\rho^\vee},
\dot w_P\cdot v^-_{\rho^\vee}\right >\\
=-\sum_{i\in I^P}\frac{1}{<w_P\cdot\rho^\vee,\alpha_i>}
(e_i^\vee)^*(u_1\cdot \rho)\left< [f_i^\vee, e_i^\vee]\dot w_P\cdot v^-_{\rho^\vee},
\dot w_P\cdot v^-_{\rho^\vee}\right >
=-\sum_{i\in I^P} (e_i^\vee)^*(u_1\cdot \rho),
\end{multline*}
using also that  $[f_i,e_i]$ acts on $\dot w_P\cdot v^-_{\rho^\vee}$ by a factor
of $<w_P\cdot \rho^\vee,\alpha_i>$.  

For the right hand summand from \eqref{e:sum} we obtain
\begin{multline*}
\left<u_1\dot w_P\dot w_0\inv u_2\inv f^\vee\cdot v^+_{\rho^\vee},\dot w_P\cdot 
v^-_{\rho^\vee}\right>\\
=\left<u_1\dot w_P\dot w_0\inv f^\vee\cdot v^+_{\rho^\vee},\dot w_P\cdot 
v^-_{\rho^\vee}\right>+
\sum_{i\in I}(e_i^\vee)^*(u_2\cdot \rho)
\left<u_1\dot w_P\dot w_0\inv [e_i^\vee, f^\vee]\cdot v^+_{\rho^\vee},\dot w_P\cdot 
v^-_{\rho^\vee}\right>\\
=\left<u_1\dot w_P\dot w_0\inv f^\vee\cdot v^+_{\rho^\vee},\dot w_P\cdot 
v^-_{\rho^\vee}\right>+
F(u_2\cdot \rho )\left<u_1\dot w_P\cdot v^-_{\rho^\vee},\dot w_P\cdot 
v^-_{\rho^\vee}\right>\\
=\sum_{i^*\in I_P}\left<u_1\dot w_P\dot w_0\inv f_i^\vee\cdot v^+_{\rho^\vee},\dot w_P\cdot 
v^-_{\rho^\vee}\right>+
F(u_2\cdot \rho ),
\end{multline*} 
by similar weight space considerations as above. Finally, using also Lemma~\ref{l:constants}~(3), the right hand summand of \eqref{e:sum} simplifies further to
\begin{multline*}
\sum_{i^*\in I_P}\left<u_1\dot w_P\dot w_0\inv f_i^\vee\cdot v^+_{\rho^\vee},\dot w_P\cdot 
v^-_{\rho^\vee}\right>+
F(u_2\cdot \rho )\\ 
= \sum_{i\in I_P}\left<u_1f_i^\vee\dot w_P\cdot v^-_{\rho^\vee},\dot w_P\cdot 
v^-_{\rho^\vee}\right>+
F(u_2\cdot \rho )\\
= -\sum_{i\in I_P}(e_i^\vee)^*(u_1\cdot \rho)\left<[e_i^\vee,f_i^\vee]\dot w_P\cdot v^-_{\rho^\vee},\dot w_P\cdot 
v^-_{\rho^\vee}\right>+
F(u_2\cdot \rho )\\
=-\sum_{i\in I_P}(e_i^\vee)^*(u_1\cdot \rho)+
F(u_2\cdot \rho ),
\end{multline*} 
noting that $<-w_P\cdot\rho^\vee,\alpha_i>=1$ for $i\in I_P$. Combining the two summands gives
\begin{multline*}
\mathcal F_P(b)= - \sum_{i\in I^P}(e_i^\vee)^*(u_1\cdot \rho)+
\left(-\sum_{i\in I_P}(e_i^\vee)^*(u_1\cdot \rho)+
F(u_2\cdot \rho )\right)\\
= F(u_2\cdot \rho)-F(u_1\cdot \rho).
\end{multline*}
\end{proof}

\begin{lem} \label{l:symmetry}
Let $Q\supseteq B$ be the parabolic subgroup determined by $W_Q=w_0 W_P w_0\inv$. 
If $b\in B^\vee_-\cap U_+^\vee (T^{\vee})^{W_P}\dot w_P\dot w_0\inv U^\vee_+$
then 
\begin{equation}
b\inv\in B^\vee_-\ \cap \ U_+^\vee (T^\vee)^{W_Q}\dot w_Q\dot w_0\inv U^\vee_+,
\end{equation} 
and the map $b\mapsto b\inv$ induces an isomorphism $\sigma_P: Z_P\to Z_Q$. Moreover we have 
\begin{equation}\label{e:symmetry}
\mathcal F_P(b)=-\mathcal F_Q(b\inv ).
\end{equation}
\end{lem}
\begin{proof}
Let us write $b=u_1 t\dot w_P\dot w_0\inv u_2\inv$ in the usual way. Then 
$$
b\inv= u_2 (\dot w_0 t\inv \dot w_0\inv)\dot w_0 \dot w_P\inv u_1\inv=
u_2 (\dot w_0t\inv \dot w_0\inv) (\dot w_0 \dot w_P\inv \dot w_0 \dot w_Q\inv) \dot w_Q \dot w_0\inv u_1\inv.
$$
Now let $\epsilon=\dot w_0 \dot w_P\inv \dot w_0\dot w_Q\inv=
\dot w_Q\inv\dot w_0\dot w_0\dot w_Q\inv$.
By Lemma~\ref{l:constants} (2) we have $ \dot w_0\dot w_0=1$, since $G^\vee$ is adjoint,
and $\epsilon=(\dot w_Q\dot w_Q)\inv\in (T^\vee)^{W_Q}$.
 
The isomorphism $\sigma_P:Z_P\to Z_Q $ is given explicitly by 
\begin{equation}\label{e:symmetrymap}
\sigma_P(t,b):=(\dot w_0 t\inv \dot w_0\inv \epsilon, b\inv).
\end{equation}  
Its inverse is $\sigma_Q$. The identity \eqref{e:symmetry} follows from Lemma~\ref{l:formula}.     
\end{proof}

Recall that by \ref{s:mirrors}.(2) we had an isomorphism
$$
B^\vee_-\cap U^\vee_+ t \dot w_P \dot w_0\inv U^\vee_+ \longrightarrow 
\mathcal R^\vee_{w_P,w_0}.
$$ 
In particular 
$Z^t_P=B^\vee_-\cap U^\vee_+ t \dot w_P \dot w_0\inv U^\vee_+$ is smooth of
dimension $n_P$. We now determine its tangent space at a point $b_0$.

\begin{lem}
Fix $t\in (T^\vee)^{W_P}$ and consider
$b_0\in B^\vee_-\cap U^\vee_+ t\dot w_P\dot w_0\inv U^\vee_+$ with 
factorization  $b_0=u_1 t\dot w_P\dot w_0\inv u_2\inv$, for $u_1,u_2\in U^\vee_+$.
We view
elements of $\mathfrak b^\vee_-$ as 
right invariant vector fields on $B^\vee_-$. Then the map
\begin{eqnarray*}
\eta:\mathfrak u^\vee_-\cap \dot w_P \cdot\mathfrak u^\vee_- &\To & \qquad \mathfrak b^\vee_-,\\
\zeta\qquad\ \  &\longmapsto &\eta_\zeta:=\pr_{\mathfrak b^\vee_- }(u_1\cdot \zeta),
\end{eqnarray*}
gives rise to an isomorphism
\begin{equation}\label{e:tangentspace}
\begin{matrix}
\mathfrak u^\vee_-\cap\dot w_P\cdot \mathfrak u^\vee_- &\To &
T_{b_0}(B^\vee_-\cap U^\vee_+ t\dot w_P\dot w_0\inv U^\vee_+),\\
\zeta\qquad \ \ &\longmapsto & \quad\qquad (\eta_\zeta)_{b_0}.
\end{matrix}
\end{equation}
\end{lem}

\begin{proof}
Let $\lambda\gg 0$ in $Q^\vee$. We consider the representations
$V(\lambda)$ and $V(\lambda+\alpha^\vee_{i^*})$ of $G^\vee$.
Then $B_-^\vee\cap U^\vee_+t\dot w_P\dot w_0\inv U^\vee_+$ inside $B_-^\vee$ is
described by the equations
\begin{eqnarray*}
\frac{\left< b\cdot v^+_{\lambda}, \dot w_P\cdot v^-_\lambda \right>}
{\left< b\cdot v^+_{\lambda+\alpha^\vee_{i^*}},
 \dot w_P\cdot v^-_{\lambda +\alpha^\vee_{i^*}}\right>}
 &=&\alpha_i^\vee(t),\\
\left< b\cdot v^+_{\lambda}, \dot w\cdot v^-_\lambda \right>\qquad&=&0 \qquad\qquad\text{for $w\in W$ with  
$w\not\ge w_P$,}
\end{eqnarray*}
where $b\in B_-^\vee$, and keeping in mind that $G^\vee$ is of adjoint type. 

Let $\zeta\in \mathfrak u_-^\vee\cap \dot w_P\cdot\mathfrak u^\vee_-$. We apply
the vector field $\eta_\zeta$ to the defining equations from above. So
\begin{multline*}
{\eta_\zeta}\left(\big\langle \variable 
\cdot v^+_\lambda,\dot w\cdot v^-_\lambda \big\rangle\right)(b_0)=
\left<\pr_{\mathfrak b^\vee_-}(u_1\cdot \zeta) b_0\cdot v^+_{\lambda}, \dot w\cdot v^-_\lambda\right>\\
=\left<u_1\zeta u_1\inv b_0\cdot v^+_{\lambda}, \dot w\cdot v^-_\lambda\right>
-\left<\pr_{\mathfrak u^\vee_+}(u_1\cdot \zeta) b_0\cdot v^+_{\lambda}, \dot w\cdot v^-_\lambda\right>\\
=\left<u_1\zeta t\dot w_P\cdot v^-_{\lambda}, \dot w\cdot v^-_\lambda\right>
-\left<\pr_{\mathfrak u^\vee_+}(u_1\cdot \zeta) u_1 t\dot w_P\cdot v^-_{\lambda}, \dot w\cdot v^-_\lambda\right>.
\end{multline*}
Since $\zeta\in\dot w_P\cdot \mathfrak u^\vee_-$ it follows that the first summand vanishes. 
The second summand is zero whenever $w\not> w_P$, by weight space considerations.
So  
$$
{\eta_\zeta}\left(\big\langle \variable 
\cdot v^+_\lambda,\dot w\cdot v^-_\lambda \big\rangle\right)(b_0)= 0
$$
if $w=w_P$ or $w\not \ge w_P$, and for any $\lambda$. In particular 
also
\begin{equation*}
{\eta_\zeta}\left(\frac{\big\langle\variable\cdot v^+_{\lambda}, \dot w_P\cdot v^-_\lambda \big\rangle}
{\big\langle\variable\cdot v^+_{\lambda+\alpha^\vee_{i^*}},
 \dot w_P\cdot v^-_{\lambda +\alpha^\vee_{i^*}}\big\rangle}\right)(b_0)=0.
\end{equation*}
It follows that $(\eta_\zeta)_{b_0}$ is tangent to $B^\vee_-\cap U^\vee_+t\dot w_P\dot w_0\inv
U^\vee_+$.

Suppose $\zeta\in\mathfrak u^\vee_-$ is homogeneous of weight $-\alpha^\vee$. Then we have 
$$
\eta_\zeta=\zeta+\bigoplus_{\beta^\vee\in\Delta^\vee_+}\mathfrak g^\vee_{-\alpha^\vee+\beta^\vee},
$$ 
and therefore $\eta$ is
injective. Comparing dimensions this implies that the map from \eqref{e:tangentspace}
is an isomorphism.   
\end{proof}

\begin{proof}[Proof of Theorem~\ref{t:main}]
Consider a fixed $b=u_1 t\dot w_P\dot w_0\inv u_2\inv$ in $B^\vee_-$, with $t\in (T^\vee)^{W_P}$ and 
$u_1,u_2\in U^\vee_+$.
\vskip .2cm
\noindent{\bf Derivatives of $\mathcal F_P$.} 
Let $\zeta\in \mathfrak u^\vee_-\cap\dot w_P\cdot \mathfrak u^\vee_-$. We may assume $\zeta$ is
homogeneous.
We want to compute the derivative of $\mathcal F_P$ in the $(\eta_\zeta)_b$ direction.  Let us write $(\eta_\zeta)_b=-\pr_{\mathfrak u^\vee_+}(u_1\cdot \zeta)+u_1\cdot \zeta$. Note that 
as for the adjoint action of $G^\vee$ on $\mathfrak g^\vee$ we also denote below 
the conjugation action of the group on itself by a dot. So $g\cdot h:=ghg\inv$ for
$g,h\in G^\vee$. Then we have
\begin{multline}\label{e:fiberder}
\eta_\zeta(\mathcal F_P)(b)
=\left .\frac{d}{ds}\right |_{0}\mathcal F_P\left ( e^{-s \pr_{\mathfrak u^\vee_+}
(u_1\cdot \zeta)}e^{s (u_1\cdot\zeta)} b\right)\\
=\left.\frac{d}{ds}\right |_{0}\mathcal F_P\left ( e^{-s \pr_{\mathfrak u^\vee_+}
(u_1\cdot \zeta)}u_1e^{s\,\zeta} t\dot w_P\dot w_0\inv u_2\inv\right)\\
=\left.\frac{d}{ds}\right |_{0}\mathcal F_P\left ( e^{-s \pr_{\mathfrak u^\vee_+}
(u_1\cdot \zeta)}u_1t\dot w_P\dot w_0\inv (\dot w_0\dot w_P\inv t\inv\cdot e^{s\zeta}) u_2\inv
\right)\\
=\left.\frac{d}{ds}\right |_{0} 
F\left(u_2\left (\dot w_0\dot w_P\inv t\inv\cdot e^{-s\zeta}\right )\cdot \rho\right)-
\left.\frac{d}{ds}\right |_{0} F\left(e^{-s \pr_{\mathfrak u^\vee_+}
(u_1\cdot \zeta)}u_1\cdot \rho\right),
\end{multline}
using Lemma~\ref{l:formula} for the last equality, and the fact that $\dot w_0\dot 
w_P\inv\cdot\zeta\in \mathfrak u^\vee_+$. The right hand summand now simplifies as follows,
\begin{multline}\label{e:right}
-\left.\frac{d}{ds}\right |_{0} F\left(e^{-s \pr_{\mathfrak u^\vee_+}
(u_1\cdot \zeta)}u_1\cdot \rho\right)=
F\left(\left [\pr_{\mathfrak u^\vee_+}(u_1\cdot \zeta), u_1\cdot \rho\right ]\right)\\
=F\left(\left [\pr_{\mathfrak u^\vee_+}(u_1\cdot \zeta),\rho\right ]\right)
=-F\left(\pr_{\mathfrak u^\vee_+}(u_1\cdot \zeta)\right)=
-F\left(u_1\cdot \zeta\right).
\end{multline}
For the left hand summand we have 
\begin{multline}\label{e:left}
\left.\frac{d}{ds}\right |_{0} 
F\left(u_2\left (\dot w_0\dot w_P\inv t\inv\cdot e^{-s\zeta}\right )\cdot \rho\right)=
-F\left(u_2\cdot \left[\dot w_0\dot w_P\inv t\inv\cdot \zeta, \rho\right ]\right)\\
=-F\left( \left[\dot w_0\dot w_P\inv t\inv\cdot \zeta, \rho\right ]\right)
=F\left( \dot w_0\dot w_P\inv t\inv\cdot \zeta\right).
\end{multline}

We now write $\mathfrak u^\vee_-=[\mathfrak u^\vee_-,\mathfrak u^\vee_-]\oplus\bigoplus_{i\in I}\mathfrak g^\vee_{-\alpha_i^\vee}$ and distinguish between two cases, corresponding to whether
$\dot w_P\inv\cdot\zeta$  lies in the one summand,
$[\mathfrak u^\vee_-,\mathfrak u^\vee_-]$, or the other,
$\bigoplus_{i\in I}\mathfrak g^\vee_{-\alpha_i^\vee}$. 
\vskip .2cm
\noindent{\bf Case 1.} Suppose $\dot w_P\inv\cdot \zeta
\in [\mathfrak u^\vee_-,\mathfrak u^\vee_-]$.  Then we have
\begin{equation*}
\dot w_0\dot w_P\inv t\inv\cdot \zeta\in\left [
\mathfrak u^\vee_+, \mathfrak u^\vee_+\right ]
\end{equation*}
and therefore $F\left(\dot w_0\dot w_P\inv t\inv\cdot \zeta\right)=0$.

\vskip .2cm
\noindent{\bf Case 2~:} In this case, since $\dot w_P\inv  \cdot\zeta$ must also lie in
$\dot w_P\inv\cdot \mathfrak u^\vee_-$, and 
$$\left(\bigoplus_{i\in I}\mathfrak g^\vee_{-\alpha_i^\vee}\right)\cap
\dot w_P\inv\cdot \mathfrak u^\vee_-=\bigoplus_{i\in I^P}\mathfrak g^\vee_{-\alpha_i^\vee},
$$
we have $\zeta\in\dot w_P\cdot \mathfrak g^\vee_{-\alpha_i^\vee}$ for some
 $i\in I^P$.
Suppose therefore $ \zeta=\dot w_P\cdot f_i^\vee$ for $i\in I^P$. 
Then 
\eqref{e:left} simplifies further to 
\begin{equation*}
F\left( \dot w_0\dot w_P\inv t\inv\cdot \zeta\right)=
F\left( \dot w_0t\inv\cdot f_i^\vee \right)=
\alpha_i^\vee(t) F\left(\dot w_0\cdot f_i^\vee \right)=
-\alpha_i^\vee(t),
\end{equation*}
using also Lemma~\ref{l:constants}~(1).

\vskip .2cm
Combining \eqref{e:right} with the above two cases for \eqref{e:left}
we get 
\begin{equation}\label{e:Fder}
\eta_\zeta(\mathcal F_P)(b)=
\begin{cases}
- F(u_1\cdot \zeta) & \text{ if $\dot w_P\inv\cdot\zeta\in [\mathfrak u^\vee_-,\mathfrak u^\vee_-] $,}\\
-\alpha_i^\vee(t)- F(u_1\dot w_P\cdot f_i^\vee)
 & \text{ if $\zeta=\dot w_P\cdot f_i^\vee$ and $i\in I_P$.}
\end{cases}
\end{equation} 
Note that, if $\eta_\zeta(\mathcal F_P)(b)=0$ for all $\zeta\in\mathfrak u^\vee_-\cap
\dot w_P\cdot\mathfrak u^\vee_-$, 
then by Case~1 above we have 
\begin{equation*}
\dot w_P\inv u_1\inv\cdot F \left |_{[\mathfrak u^\vee_-,\mathfrak u^\vee_-]}\right .=0.
\end{equation*}
This implies that $u_1\dot w_P B^-/B^-$, or equivalently $b\dot w_0 B^-/B^-$, 
lies in the non-equivariant
Peterson variety $Y_P$. 
\vskip .2cm

\noindent{\bf Logarithmic derivative of $\phi$.} Let $\zeta\in \mathfrak u^\vee_-\cap\dot w_P\cdot \mathfrak u^\vee_-$. Decomposing
\begin{equation*}
(\eta_{\zeta})_b=\pr_{\mathfrak u^\vee_-}(u_1\cdot \zeta)+
\pr_{\mathfrak h^\vee} (u_1\cdot \zeta)
\end{equation*}
we see that $(\eta_\zeta)_b$ lifts to the tangent vector
$$
(\tilde\eta_\zeta)_{(b,h_R^\vee)}=(\pr_{\mathfrak b^\vee_-}(u_1\cdot \zeta),
\pr_{\mathfrak h^\vee} (u_1\cdot \zeta))
$$
in $T_{(b,h^\vee_R)}(\tilde Z_P^t)\subset \mathfrak b^\vee_-\oplus\mathfrak h^\vee$. The 
logarithmic derivative of $\tilde \phi$ in this direction is therefore given by
\begin{multline}\label{e:phiderivative}
\tilde\eta_{\zeta}(\ln \tilde\phi(\ ;h))(b,h^\vee_R)=\left.\frac{d}{ds}\right |_{s=0}<h, s\pr_{\mathfrak h^\vee}(u_1\cdot \zeta)+h^\vee_R>\\
=<h,\pr_{\mathfrak h^\vee}(u_1\cdot \zeta)>.
\end{multline}
Note again that $\tilde\eta_{\zeta}(\ln \tilde\phi(\ ;h))(b,h^\vee_R)$ no longer depends on the
choice of lift $(b,h^\vee_R)\in \tilde Z^t_P$ of $b$, and is a well-defined function on $Z^t_P$. 
We view this as the derivative of the multi-valued
function $\ln\phi (\ ;h)$ at the point $b$ in $Z^t_P$ in the direction $(\eta_\zeta)_b$, and 
may also denote it by $\eta_\zeta(\ln\phi)(b;h)$.

\vskip .2cm
\noindent{\bf The critical points  of $\mathcal F_P+\ln\phi (\ ;h)$ along fibers.}
By definition
\begin{equation*}
Z^{crit}_{P,T}=\left\{(t,b;h)\in Z_P\x\mathfrak h\ |\ 
\eta_{\zeta}(\mathcal F_P)(b)+\eta_{\zeta}(\ln\phi)(b;h)=0\ \text{ for all 
$\zeta\in\mathfrak u^\vee_-\cap \dot w_P\cdot \mathfrak u^\vee_-$}\right\}.
\end{equation*}
As before we have two cases for $\zeta$.
\begin{enumerate}
\item
If
$\zeta\in \dot w_P\cdot [\mathfrak u^\vee_-,\mathfrak u^\vee_-]
\cap \mathfrak u^\vee_-$
then, by \eqref{e:phiderivative} and \eqref{e:Fder},  
 \begin{equation*}
\eta_\zeta(\mathcal F_P)(b) + \eta_\zeta(\ln\phi)(b;h)= u_1\inv\cdot (-F+h)\ (\zeta),
 \end{equation*} 
where $h\in\mathfrak h=(\mathfrak h^\vee)^*$ is considered as 
an element of $(\mathfrak g^\vee)^*$. Let us replace $\zeta$ by 
$\bar \zeta :=\dot w_P\inv\cdot  \zeta$, so  $\bar
 \zeta\in 
[\mathfrak u^\vee_-,\mathfrak u^\vee_-]\cap \dot w_P\inv\cdot\mathfrak u^\vee_-$. 
The critical point condition $\eta_\zeta(\mathcal F_P+\ln\phi (\ ;h))(b)=0$ 
in this case reads 
\begin{equation*}
  \dot w_P\inv u_1\inv\cdot (-F+h)(\bar \zeta)=0, \quad\text{ for all $\bar \zeta\in
[\mathfrak u^\vee_-,\mathfrak u^\vee_-]\cap \dot w_P\inv\cdot \mathfrak u^\vee_-$,}
\end{equation*} 
in terms of $\bar \zeta$.
\item
If $\zeta= \dot w_P\cdot f_i^\vee$ for $i\in I^P$, then by \eqref{e:phiderivative}
and \eqref{e:Fder} we have 
\begin{equation*}
\eta_{\dot w_P\cdot f_i^\vee}(\mathcal F_P)(b)+\eta_{\dot w_P\cdot f_i^\vee}(\ln\phi)(b;h)
 =
-\alpha_i^\vee(t) - (u_1\inv\cdot (F - h))
(\dot w_P\cdot f_i^\vee).
\end{equation*}
\end{enumerate}

Finally, note that if $\bar\zeta\in [\mathfrak u^\vee_-,\mathfrak u^\vee_-]\cap \dot w_P\inv\cdot
\mathfrak u^\vee_+$, then $\dot w_P\cdot\bar\zeta\in 
[\mathfrak u^\vee_+,\mathfrak u^\vee_+]$, so 
$$
\dot w_P\inv u_1\inv\cdot (-F+h)(\bar\zeta)=0
$$ 
automatically. Therefore the critical point condition (1) implies 
$$
\dot w_P\inv u_1\inv\cdot (-F+h)|_{[\mathfrak u^\vee_-,\mathfrak u^\vee_-]}=0.
$$ 

Combining (1) and (2) above, we find that the critical point locus $Z^{crit}_{P,T}$ is given by
\begin{equation}\label{e:critlocus}
Z^{crit}_{P,T}=\left\{ (t,b;h)\in Z_P\x\mathfrak h\ \left |
\ \begin{matrix}\dot w_P\inv u_1\inv\cdot (F-h)\ |_{[\mathfrak u^\vee_-,\mathfrak u^\vee_-]}=0
\\ \text{and}
\\
(F-h)(u_1\dot w_P\cdot f^\vee_i)=-\alpha_i^\vee(t), \\
\text{for $b=u_1t \dot w_P\dot w_0\inv u_2\inv$ and $i\in I^P$}
\end{matrix}
\right.
\right\}.
\end{equation}
This implies that $b\mapsto b\dot w_0 B^\vee_-/B^\vee_-=u_1\dot w_P\cdot B^\vee_-/B^\vee_-$ 
defines a map
\begin{equation}\label{e:desiredmap}
Z^{crit}_{P,T}\To \mathcal Y_P^*.
\end{equation}
Comparing with Peterson's description of the quantum parameters, Section~\ref{s:qparams}, 
we see that the diagram~\eqref{e:diagram} commutes. 

To show that  \eqref{e:desiredmap} is an isomorphism, consider
$(u\dot w_PB^\vee_-/B^\vee_-,h)\in \mathcal Y_P^*$, 
where $u\in U_-^\vee$.  Define $t\in (T^\vee)^{W_P}$ by the condition
\begin{equation*}
(F-h)(u\dot w_P\cdot f^\vee_i)=-\alpha_i^\vee(t),\qquad \text{for all $i\in I^P$.}
\end{equation*}
Then since $u\dot w_PB^\vee_-/B^\vee_-\in \mathcal R^\vee_{w_P,w_0}$ there is a
unique $b=u_1t\dot w_P\dot w_0\inv u_2\inv\in Z_P^t$ with 
$$
b\dot w_0B^\vee_-/B^\vee_-=u\dot w_P B^\vee_-/B^\vee_-,
$$ 
as in (2) of Section~\ref{s:mirrors}. It is clear from \eqref{e:critlocus} that $(t,b;h)\in Z^{crit}_{P,T}$ and so 
we have defined an inverse to \eqref{e:desiredmap}.

The description
\eqref{e:stabilizer} of $Z^{crit}_{P,T}$ is an immediate consequence of \eqref{e:desiredmap} being an isomorphism,
together with the analogous result for 
the Peterson variety from Section~\ref{s:stabilizer}, which is due to Dale Peterson
and originates from a description of Kostant's for the leaves of 
the Toda lattice. 
We explain the proof here for completeness. It starts with the observation that 
the condition
$$
\dot w_P\inv u_1\inv\cdot (F-h)\ |_{[\mathfrak u^\vee_-,\mathfrak u^\vee_-]}=0
$$
implies that the element
$$
b\inv\cdot (F-h)= u_2\dot w_0\dot w_P\inv t \inv u_1\inv\cdot (F-h)
$$
vanishes on $[\mathfrak u^\vee_+,\mathfrak u^\vee_+]$, and therefore
lies in $\bigoplus_{i\in I}(\mathfrak g^\vee_{\alpha_i})^*\oplus (\mathfrak h^\vee)^*$.
Now the fact that $b\in B^\vee_-$ implies immediately that 
$$
b\inv\cdot (F-h)\ (h^\vee)=(F-h)\ (b \cdot h^\vee)=-<h,h^\vee>
$$
for any $h^\vee\in\mathfrak h^\vee$.
And direct calculation 
using  the second condition from \eqref{e:critlocus}, along with the fact that $\alpha_i(t)=1$
if $i\in I_P$, shows that
$$
b\inv\cdot (F-h)\ (e_i^\vee)=1
$$
for any $i\in I$. Therefore in fact
$b\inv\cdot (F-h)=F-h$ and 
\begin{equation}\label{e:stabilizerinclusion}
Z^{crit}_{P,T} \subseteq \{ (t,b;h)\in Z_P\ |\ b\cdot(F-h)=F-h\}.
\end{equation}
The opposite inclusion follows using the identity
$$
t \dot w_0\inv u_2\inv \cdot (F-h)=t \dot w_0\inv u_2\inv b\inv\cdot (F-h)=\dot w_P\inv u_1\inv\cdot (F-h),
$$
for any $b$ in the right hand side of  \eqref{e:stabilizerinclusion}.
\end{proof}

\section{Deodhar stratifications and standard coordinates}\label{s:Deodhar}
In this section we introduce coordinate systems on intersections of
opposite Bruhat cells. These will be used in the subsequent sections firstly 
to define the holomorphic $n_P$-forms
we need to state the mirror conjecture for $G/P$, and secondly 
to compare our mirror construction with the one from 
\cite{Giv:MSFlag,JoeKim:EquivMirrors} in type $A$.

\subsection{}
Intersections of opposed Bruhat cells 
$\mathcal R_{v,w}$ were decomposed  
into strata isomorphic to products of the form 
$\C^p\x(\C^*)^q$ by Deodhar \cite{Deo:decomp}. 
We will give a 
practical definition of these strata following 
\cite{MarRie:ansatz}. This latter description
has the advantage
of providing for every stratum natural 
coordinates to work with. 

Let $w\in W$ and $s_{i_1} s_{i_2} \dotsc s_{i_m}=w$
be a fixed reduced expression 
which we denote by $\mathbf i=(i_1,\dotsc,i_m)$.   
We consider a sequence of integers 
$1\le j_1<\dotsc< j_t\le m$ as 
giving a subexpression 
$s_{i_{j_1}}\dotsc s_{i_{j_t}}$ of $s_{i_1}\dotsc s_{i_m}$.
We say it is a subexpression for 
$v$ if $s_{i_{j_1}}\dotsc s_{i_{j_t}}=v$. 
Note that $(i_{j_1},\dotsc, i_{j_t})$ need not be
a reduced expression of $v$.

A subexpression $\mathbf j=(j_1,\dotsc, j_t)$ of $\mathbf i$ 
is called {\it distinguished} if 
 $$(s_{i_{j_1}}\dotsc s_{i_{j_l}}) s_{i_k} > s_{i_{j_1}}\dotsc s_{i_{j_l}}
 \ \text{for all}\ j_l< k< j_{l+1},$$ 
where $1\le l\le t$. There is a unique subexpression 
for $v$ with the stronger property that 
 $$( s_{i_{j_1}}\dotsc s_{i_{j_l}} ) s_{i_k} >  s_{i_{j_1}}\dotsc s_{i_{j_l}}
\ \text{for all $j_l<k \le j_{l+1}$,}
$$ 
where $1\le l\le t$. We may set $j_{t+1}=m+1$ everywhere above. 
We call this subexpression the {\it positive
subexpression} for $v$. It is the unique distinguished 
subexpression that gives a reduced expression for $v$. 

Deodhar's construction
associates to any reduced expression of
$w$ a stratification of $\mathcal R_{v,w}$ which has
a stratum for every distinguished 
subexpression for $v$. And the positive subexpression for $v$ 
corresponds to the unique open stratum.

For a reduced expression $\mathbf i$ and subexpression 
$\mathbf j$ let
 \begin{eqnarray*}
J_0(\mathbf j)&=&\{1,\dotsc,m\}\setminus\{j_1,\dotsc, j_t\},\\
J_+(\mathbf j)&=&\left\{l \ \left | \ l=j_r \text{ some $r=1,\dotsc, t$, and } 
 s_{i_{j_1}}\dotsc s_{i_{j_{r}}}>
s_{i_{j_1}}\dotsc s_{i_{j_{r-1}}} \right .\right\},\\
J_-(\mathbf j)&=&\left\{l \ \left | \ l=j_r \text{ some $r=1,\dotsc, t$, and } 
 s_{i_{j_1}}\dotsc s_{i_{j_{r}}}<
s_{i_{j_1}}\dotsc s_{i_{j_{r-1}}} \right .\right\},
\end{eqnarray*}
where we suppress the $\mathbf i$ in the notation since it is
usually clear from context.  If $\mathbf j$ is distinguished, 
define a subset  
$\mathcal R_{\mathbf j,\mathbf i}$ of the flag variety
$G/B_-$ by
$$\mathcal R_{\mathbf j,\mathbf i}:=
\left\{g_1 \dotsc g_m B_-/B_- \ \left |\ 
g_l=\begin{cases}
x_{i_l}(t_{l}),\ 
\ t_l\in\C^*,& \text{if}\ l\in J_0(\mathbf j) \\
\dot s_{i_l},   & 
\text{if}\ l\in J_+(\mathbf j)\\ 
y_{i_l}({m_l})\dot s_{i_l}\inv,\ m_l\in \C,& \text{if}\ 
l\in J_-(\mathbf j)
\end{cases} \right.
\right\}.
$$
Here the parameters $t_l\in\mathbb C^*$ and 
$m_l\in\mathbb C$ can also be used as coordinates
on $\mathcal R_{\mathbf j,\mathbf i}$ giving an 
isomorphism $\mathcal R_{\mathbf j,\mathbf i}
\overset\sim\longrightarrow 
(\mathbb C^*)^{|J_0(\mathbf j) |}\times
\mathbb C^{|J_-(\mathbf j)|}$. We will refer
to these coordinates as the standard coordinates
on $\mathcal R_{\mathbf j,\mathbf i}$.
If $\mathbf j_+$ is the positive subexpression 
for $v$ in $\mathbf i$
then 
$\mathcal R_{\mathbf j_+,\mathbf i}\overset\sim\To 
(\mathbb C^*)^{\ell(w)-\ell(v)}$.  

By \cite[Proposition~5.2]{MarRie:ansatz} 
the $\mathcal R_{\mathbf j,\mathbf i}$ agree   
precisely with Deodhar's strata in $\mathcal R_{v,w}$.
So fixing $\mathbf i$  we have 
 $$
\mathcal 
R_{v,w}=
\bigsqcup_{\mathbf j}\mathcal 
R_{\mathbf j,\mathbf i},
$$
where the union is over all distinguished subexpressions $\mathbf j$
of $\mathbf i$.
Note only that our conventions differ from \cite{MarRie:ansatz}
in that $B_+$ and $B_-$ are interchanged.

\section{A holomorphic $n_P$-form on $\mathcal R^\vee_{w_P,w_0}$}

To define the oscillatory 
integrals and state the mirror conjecture for $G/P$
we require holomorphic $n_P$-forms 
on the fibers of the proposed mirror family. 
Therefore we want to define a holomorphic $n_P$-form on 
an intersection of opposed Bruhat cells $\mathcal R^\vee_{w_P,w_0}$. 
This holomorphic differential form will be defined by writing it down
explicitly on a large enough open subset of $\mathcal R^\vee_{w_P,w_0}$.
\begin{footnote}
{M.~Brion pointed out to us that by \cite{Brion:Lectures} Theorem~4.2.1(i) 
the canonical bundle of the closure of any $\mathcal R_{v,w}$, the so-called
 Richardson variety $X_{v,w}:=\overline{\mathcal R_{v,w}}$, is 
 $\mathcal O_{X_{v,w}}(\partial X_{v,w})$, and he conjectured that our form
 might come from there by restriction, \cite{Brion:June07}. If so, this would
 give a more intrinsic definition of our form, at least up 
 to scalar. } 
\end{footnote}

Let $\mathbf i$ be
a reduced expression of $w_0$ and
$\mathbf j=\mathbf j_+(\mathbf i)$ the corresponding
positive subexpression for $w_P$. Consider the 
open Deodhar stratum 
 \begin{equation}\label{e:openpart}
\mathcal R^\vee_{\mathbf j_+(\mathbf i),\mathbf i}=
\left\{g_1\dotsc g_N B^\vee_-/B^\vee_- \left | g_l=
\begin{cases} x_{i_l}(t_l)\ \text{for}\ t_l\in\C^* & 
\text{if}\ l\in J_0(\mathbf  j_+(\mathbf i)) \\
\dot s_{i_l} & \text{otherwise.}
\end{cases}\right .
\right\},
 \end{equation}
in $\mathcal R^\vee_{w_P,w_0}$. Let $\mathcal U$ be the union
of these open sets. So 
 $$
\mathcal U=\bigcup_{\mathbf i} \mathcal R^\vee_{\mathbf j_+(\mathbf i),\mathbf i},
$$
where $\mathbf i$ ranges over all the reduced expressions of $w_0$.

\begin{lem}[essentially Lemma~3.6 in \cite{Zel:ConnComps}]\label{l:U}
$\mathcal U$ is an open dense subset of 
$\mathcal R^\vee_{w_P,w_0}$ with complement of codimension greater than or equal 
to $2$.
\end{lem}
 
\begin{proof} Since $\mathcal R^\vee_{w_P,w_0}$ is irreducible
it is clear that $\mathcal U$ is open dense. 
We have an isomorphism 
 \begin{align*}
U^\vee_-\cap B^\vee_+ \dot w_P \dot w_0 B^\vee_+ \overset \sim\To 
\mathcal R^\vee_{w_P,w_0}\ :\quad
 u\mapsto u\dot w_0 B^\vee_-,
 \end{align*} 
whereby the double Bruhat cell $B^\vee_-\cap B^\vee_+ \dot w_P\dot w_0 B^\vee_+$
in the group can be identified with $\mathcal R^\vee_{w_P,w_0}\x T^\vee$. 
Now Lemma~3.6 in \cite {Zel:ConnComps}, 
which is about $B^\vee_-\cap B^\vee_+\dot w_P\dot w_0
B^\vee_+$, implies the lemma.
\end{proof}

\begin{prop}\label{p:omega} Fix a reduced expression $\mathbf i_0$ of $w_0$. 
There is a unique holomorphic $n_P$-form $\omega$ on 
$\mathcal R^\vee_{w_P,w_0}$ 
such that the 
restriction of $\omega$ to $\mathcal R^\vee_{\mathbf j_+(\mathbf i),\mathbf i}$
is given by
 $$
\omega|_{\mathcal R^\vee_{\mathbf j_+(\mathbf i),\mathbf i}}=\epsilon_{\mathbf i}\ 
\bigwedge_{l\in J_0(\mathbf j_+(\mathbf i))} \frac{d t_l}{t_l},
 $$ 
in terms of the 
standard coordinates $t_l$ on 
$\mathcal R^\vee_{\mathbf j_+(\mathbf i),\mathbf i}$, where
$\epsilon_{\mathbf i}\in \{\pm 1\}$ and $\epsilon_{\mathbf i_0}=1$. 
Here we use the obvious order on $J_0(\mathbf j_+(\mathbf i))$
for defining the wedge product. 
\end{prop}

\begin{proof} By Lemma~\ref{l:U} and Hartog's theorem if $\omega$ is
well defined on $\mathcal U$ then it extends holomorphically to 
all of $\mathcal R^\vee_{w_P,w_0}$.

Let $\mathbf i$ and $\mathbf i'$ be reduced expressions 
of $w_0$ such that $\mathbf i$ is obtained from $\mathbf i'$
by a single braid relation. 
It suffices to show that the rational transformation 
$(t_1,\dotsc, t_{n_P})\mapsto(t'_1,\dotsc,t'_{n_P})$
from the standard coordinates on  
$\mathcal R^\vee_{\mathbf j_+(\mathbf i),\mathbf i}$ to those of
$\mathcal R^\vee_{\mathbf j_+(\mathbf i'),\mathbf i'}$ gives
 $$
\frac{dt_1}{t_1}\wedge\cdots\wedge\frac{dt_{n_P}}{t_{n_P}}
=\pm\frac{dt'_1}{t'_1}\wedge\cdots\wedge\frac{dt'_{n_P}}{t'_{n_P}}.
 $$ 

The remainder of the proof consists of checking the
possible coordinate transformations that can occur. 
\vskip .2cm
\noindent{\bf Simply laced case~:}
\vskip .1cm
\begin{enumerate} 
\item 
If $s_i s_j= s_j s_i$ then 
\begin{equation*}
x_i(a) x_j(b)=x_j(b) x_i(a)
\end{equation*}
giving the simplest change of coordinates $C_0(a,b)=(b,a)$. 
\item
If $s_i s_j s_i=s_j s_i s_j$ then it is easy to check that
\begin{equation*}
x_i(a) x_j(b) x_i(c)= x_j\left(\frac{bc}{a+c}\right) x_i(a+c) 
x_j\left(\frac{ab}{a+c}\right)
\end{equation*}
and 
\begin{equation*}
x_i(a) x_j(b)\dot s_i= x_j(b )\dot s_i x_j(a b) y_i(-a).
\end{equation*}
We may record these two changes of coordinates as 
\begin{eqnarray*}
C_1(a,b,c)&=&\left(\frac{bc}{a+c}, a+c, \frac{ab}{a+c}\right)\\
C_2(a,b)&=&(b, a b)
\end{eqnarray*}
Note also that
\begin{equation*}
x_i(a)\dot s_j\dot s_i=\dot s_j\dot s_i x_j(a).
\end{equation*}
\end{enumerate}
\vskip .1cm
\noindent{\bf Type $B_2$ braid relations~:}
\vskip .1cm

If $s_i s_j s_i s_j=s_j s_i s_j s_i$ and $\al_i$ is the long root, 
then the following
relation holds (see \cite[Section 3.1]{BeZel:TotPos}).
\begin{align*} &\quad x_j(a) x_i(b) x_j(c) x_i(d)=x_i(d') x_j(c') x_i(b')
x_j(a')\\ &{\text\rm where}\quad a'=\frac {a b c}y \qquad b'=\frac
{y^2}x \qquad c'=\frac xy  \qquad  d'=\frac{b c^2
d}{x}\qquad\qquad\\ &{\text\rm and}\qquad x= a^2 b + d(a+c)^2
\qquad y= a b + d (a+c).
\end{align*}
Let us denote this change of coordinates by
\begin{equation*}
C_3(a,b,c,d)=(d',c',b',a').
\end{equation*}
Its inverse is $(C_3)^{-1}(d,c,b,a)=(a',b',c',d')$, where
$a',b',c',d'$ are given by the same
formulas as above. 

Furthermore it is easy to check the pairs of (inverse) identities 
\begin{eqnarray*}
x_i(a) x_j(b) x_i(c)\dot s_j &= &x_j\left(\frac{b c}{a+c}\right) 
x_i(a+c)\dot s_j
x_i\left(\frac{a b^2 c}{a+c}\right) y_j\left(\frac{-a b}{a+c}\right),\\
x_j(a) x_i(b) \dot s_j x_i(c) &= &
x_i\left(\frac{b c}{c + a^2 b}\right) x_j\left(\frac{c+a^2 b}{a b}\right) 
x_i\left(\frac{a^2 b^2}{c+ a^2 b}\right)\dot s_j y_i\left(\frac{c}{a b}\right),
\end{eqnarray*}
and
\begin{eqnarray*}
x_j(a) x_i(b) \dot s_j \dot s_i
&=&x_i(b)\dot s_j \dot s_i x_j(a b) u,\\
x_i(a) \dot s_j \dot s_i x_j(b)
&=& x_j\left(\frac ba\right) x_i(a)\dot s_j\dot s_i u',
\end{eqnarray*}
for some $u$ and $u'$ in $U_-$, and 
\begin{eqnarray*}
x_j(a)x_i(b) x_j(c)\dot s_i= x_i\left(\frac{ c^2 b}{(a+c)^2}\right)
x_j(a+c)\dot s_i x_j\left(\frac{a b c}{a+c}\right)
y_i\left(\frac{-a^2b-2 a b c}{(a+c)^2}\right)\\
x_i(a) x_j(b)\dot s_i x_j(c)=x_j\left (\frac{b c}{c+a b}\right)
x_i\left(\frac{(c+a b)^2}{a b^2}\right) x_j\left(\frac{a b^2}{c+ a b}
\right)\dot s_i y_i\left(\frac{2 a b c + c^2}{a b^2}\right),
\end{eqnarray*}
and
\begin{eqnarray*}
x_i(a) x_j(b) \dot s_i \dot s_j&=&x_j(b)\dot s_i \dot s_j x_i(a b^2) y,\\
x_j(a) \dot s_i \dot s_j x_i(b)&=
& x_i\left(\frac b{a^2}\right) x_j(a)\dot s_i\dot s_j y'
\end{eqnarray*}
for some $y$ and $y'$ in $U_-$. 
Finally 
\begin{eqnarray*}
x_i(a) \dot s_j\dot s_i\dot s_j&=&\dot s_j \dot s_i\dot s_j x_i(a),\\
x_j(a) \dot s_i\dot s_j\dot s_i&=&\dot s_i \dot s_j\dot s_i x_j(a).
\end{eqnarray*}
We record the remaining four nontrivial changes of coordinates 
\begin{eqnarray*}
C_4(a,b,c)\ &=&\left(\frac{b c}{a+c}, a+c,  \frac{a b^2 c}{a+c}\right),\\
C_5(a,b)\ \ &=& \left(b, a b\right),\\
C_6(a,b,c)\ &=&\left(\frac{c^2 b}{(a+c)^2}, a+c, \frac{a b c}{a+c}\right),\\
C_7(a,b)\ \ &=& \left(b, a b^2\right).
\end{eqnarray*} 

\vskip .2cm
\noindent{\bf Type $G_2$ braid relations~:}
\vskip .2cm
\begin{enumerate}
\item
 If $s_i s_j s_i s_j s_i s_j=s_j s_i s_j s_i s_j s_i$ and
$\al_i$ is the long root, then we have (see
\cite[Section 3.1]{BeZel:TotPos})
$$
x_i(a) x_j(b) x_i(c) x_j(d) x_i(e) x_j(f)= x_j(f') x_i(e') x_j(d')
x_i(c') x_j(b') x_i(a')
$$
for
$$
\begin{aligned}
&a'=\frac{a b c^2 d e}{\pi_1} \qquad &b'=\frac{\pi_1^3}{\pi_4}
\qquad &c'=\frac{\pi_4}{\pi_1 \pi_2} \qquad\qquad\qquad\qquad \\
&d'=\frac{\pi_2^3}{\pi_3 \pi_4} \qquad &e'=\frac{\pi_3}{\pi_2}
\qquad &f'=\frac{b c^3 d^2 e^3 f}{\pi_3}
\end{aligned}
$$
Where
\begin{align*}
\pi_1= &\ abc^2d+ab(c+e)^2 f+ (a+c)de^2f  \\
\pi_2= &\ a^2 b^2 c^3 d + a^2 b^2(c+e)^3 f + (a+c)^2 d^2 e^3 f
+ abde^2 f (3ac+2c^2+2ce+2ae) \\
\pi_3= &\ a^3 b^2 c^3 d + a^3 b^2(c+e)^3 f + (a+c)^3 d^2 e^3 f
+ a^2 bde^2 f (3ac+3c^2+3ce+2ae) \\
\pi_4= &\ a^2 b^2 c^3 d(abc^3 d + 2ab(c+e)^3 f +
(3ac+3c^2+3ce+2ae)de^2 f) +\ \\
       &\ f^2(ab(c+e)^2+(a+c)d e^2)^3.
\end{align*}
This gives rise to the eighth relevant 
change of coordinates 
$$C_8(a,b,c,d,e,f)=(f',e',d',c',b',a').$$ 
\item Next we have the relation
\begin{multline*}
x_i(a)x_j(b) x_i(c) x_j(d) \dot s_i x_j(e)=\\
=x_j\left(\frac{b c^3 d^2 e}{x}\right) x_i\left(\frac xy\right) 
x_j\left(\frac{y^3}{x z}\right) 
x_i\left(\frac{z}{a b c^2 d v}\right) 
x_j\left(\frac{a^3 b^3 c^6 d^3}{z}\right) \dot s_i u,
\end{multline*}
where $u\in U_-$ and
\begin{align*}
x&= 3 a c d e(c d + a b + a d)+ c^3 d^2 e+ a^3(2 b d e+ d^2 e+ b^2(c^3 d+e)),\\
y&= c d e(c d+a b+ a d)+ a c d(b+d) e+ a^2(2 b d e+ d^2+b^2(c^3 d+e)),\\
z&=e(c d + a b+ a d)(y+ a^2 b^2 c^3 d)+ a^2 b^2 c^4 d^2(e+a b c^2),\\
v&=a^2 b^2 c^3 d+ e(c d+ a b+ a d)^2.
\end{align*}
So we may take
$$
C_9(a,b,c,d,e)=\left(\frac{b c^3 d^2 e}{x},\frac xy,\frac{y^3}{x z},
\frac{z}{a b c^2 d v},\frac{a^3 b^3 c^6 d^3}{z}\right).
$$

Similarly,
 \begin{multline*}
x_j(a)x_i(b) x_j(c) x_i(d) \dot s_j x_i(e)\\=
x_i\left(\frac{b c d^2 e}{x'}\right) x_j\left(\frac{{x'}^3}{y'}\right) x_i\left(\frac{y'}{x' z'}\right) 
x_j\left(\frac{{z'}^3}{a b^3 c^3 d^3 y'}\right) 
x_i\left(\frac{a b^3 c^2 d^3}{z'}\right) \dot s_j u
 \end{multline*}
where $u\in U_-$ and
\begin{align*}
x'&= c d^2 e+ a(2 b d e+ d^2 e+ b^2(c d+e)),\\
y'&= c^2 d^6 e^3+a^2(2 b d e+ d^2 e+ b^2(c d+e))^3+ a c d^3(b+d)e^2
(4 b d e+ 2 d^2 e+ b^2(3 c d+ 2 e)),\\
z'&=c d^3 e^2+ a(3 b d^2 e^2+ d^3 e^2+b^3(c d+e)^2+ b^2 d e (2 c d+3 e)).
\end{align*}
This gives
 $$
C_{10}(a,b,c,d,e)=
\left(\frac{b c d^2 e}{x'},\frac{{x'}^3}{y'},\frac{y'}{x' z'},
\frac{{z'}^3}{a b^3 c^3 d^3 y'},\frac{a b^3 c^2 d^3}{z'}\right).
 $$

\item Next we have
 \begin{multline*}
x_j(a)x_i(b)x_j(c)x_i(d)\dot s_j \dot s_i=\\
x_i\left(\frac{b c^3 d^2}{a z_1^2+ c d z_2}\right)
x_j\left(\frac{ a z_1^2 + c d z_2}{z_1^2}\right)
x_i\left(\frac{z_1^3}{a z_1^2+ c d z_2}\right)
\dot s_j\dot s_i x_j\left(\frac{a b c^2 d}{z_1}\right) u
 \end{multline*}
where
$z_1=c d + a(b+d),\ 
z_2=a^2 b+(a+c)^2d$ 
and $u\in U_-$,
and also
\begin{multline*}
x_j(a)x_i(b)x_j(c)\dot s_i \dot s_j x_i(d)=\\
x_i\left(\frac{b c^3 d}{z_4}\right)
x_j\left(\frac{z_4}
{z_3}\right)
\dot s_i \dot s_j
x_i\left(\frac{z_3^3}{a^3 b^3 c^6 z_4}\right)
x_j\left(\frac{a^2 b^2 c^4}{z_3}\right) u'
 \end{multline*}
where 
$$z_3= a^2 b^2 c^3+(a+c)^2 d\quad,
z_4=a^3 b^2 c^3+ (a+c)^3 d,$$
and $u'\in U^-$.
So we set 
 \begin{align*}
C_{11}(a,b,c,d)&=\left(\frac{b c^3 d^2}{a z_1^2+ c d z_2},
\frac{ a z_1^2 + c d z_2}{z_1^2},
\frac{z_1^3}{a z_1^2+ c d z_2},
\frac{a b c^2 d}{z_1}\right),\\
C_{12}(a,b,c,d)&=\left(\frac{b c^3 d}{z_4},
\frac{z_4}{z_3},\frac{z_3^3}{a^3 b^3 c^6 z_4},
\frac{a^2 b^2 c^4}{z_3}
\right).
 \end{align*}
\item
Moreover
\begin{equation*}
x_j(a)x_i(b)x_j(c)\dot s_i\dot s_j\dot s_i=
x_i\left(\frac{b c^3}{(a+c)^3} \right)
x_j\left(a+c\right)
\dot s_i \dot s_j \dot s_i
x_j\left(\frac{a b c^2}{a+c}\right) u
\end{equation*}
and similarly
\begin{equation*}
x_j(a)x_i(b)\dot s_j\dot s_i\dot s_jx_i(c)=
x_i\left(\frac{b c}{a^3 b^2+c} \right)
x_j\left(\frac{a^3 b^2+c}{a^2 b^2}\right)
x_i\left(\frac{a^3 b^3}{a^3 b^2+c}\right)
\dot s_j \dot s_i \dot s_j u'
\end{equation*}
for some $u,u'\in U^-$.
So
\begin{eqnarray*}
C_{13}(a,b,c)&=\left(\frac{b c^3}{(a+c)^3}
,a+c,\frac{a b c^2}{a+c}\right),\\
C_{14}(a,b,c)&=\left(\frac{b c}{a^3 b^2+c},\frac{a^3 b^2 +c}{a^2b^2},
\frac{a^3 b^3}{a^3 b^2+c}\right).
\end{eqnarray*}
\item
Finally,
\begin{equation*}
x_i\left(a\right)
x_j\left(b\right)
\dot s_i
\dot s_j \dot s_i \dot s_j =
x_j(b)\dot s_i\dot s_j\dot s_i\dot s_jx_i(a b^3)u
\end{equation*}
and
\begin{equation*}
x_j(a)x_i(b)\dot s_j\dot s_i\dot s_j\dot s_i=
x_i\left(b \right)
\dot s_j
\dot s_i
\dot s_j \dot s_i x_j(a b) u'
\end{equation*}
up to $u,u'\in U^\vee_-$. 
Also we have
\begin{eqnarray*}
x_i(a)\dot s_j\dot s_i\dot s_j\dot s_i\dot s_j&=
\dot s_j\dot s_i\dot s_j\dot s_i\dot s_j x_i(a)\\
\dot s_i\dot s_j\dot s_i\dot s_j\dot s_i x_j(a)&=
x_j(a)\dot s_i\dot s_j\dot s_i\dot s_j \dot s_i
\end{eqnarray*}
So the last new coordinate transformation is
\begin{eqnarray*}
C_{15}(a,b)&=\left(b,ab^3\right).
\end{eqnarray*}
\end{enumerate}
\vskip .2cm
Here the transformations in (2-5) immediately
above 
were computed with the help of Mathematica, 
realizing $G_2$ inside a group of type $B_3$
and using all of the relations from types 
$A_2$ and $B_2$.

Now let
$$
L:\mathbf t=(t_1,\dotsc, t_m)\mapsto (L^1(\mathbf t),\dotsc, L^m(\mathbf t))
$$ 
be one of the 
changes of coordinates $C_j$ with $j=0,\dotsc, 15$. The form given by
$$
\bigwedge_{i=1}^l \frac{dt_i}{t_i}
$$
is invariant up to sign under 
these changes of coordinates if for each of the $L=C_j$
 $$
\frac{Jac(L)}{L^1\cdot\dotsc\cdot L^m}=\pm\frac{1}{t_1\cdot\dotsc\cdot t_m},
$$ 
where 
$Jac(L)= \det\left(\frac{\partial L^i}{\partial t_k}\right)_{i,k=1,\dotsc, m}$
is the Jacobian. 
This is the case as can easily be checked e.g. using Mathematica.
The sign is minus precisely in the cases
$$C_0, C_2, C_3, C_5, C_7, C_8, C_{11}, C_{12}, C_{15},$$
where there is an even number of coordinates involved in the coordinate
transformation.
\end{proof}

\begin{rem}
Notice that all of the coordinate transformations $C_0,\dotsc, C_{15}$
are subtraction-free rational functions. The 
well-defined subset in the real
points of an intersection of Bruhat cells  $\mathcal R_{v,w}$, consisting of those points in an (any) 
open Deodhar stratum 
$\mathcal R_{\mathbf j_+(\mathbf i),\mathbf i}$ all of whose 
canonical coordinates 
take values in $\mathbb R_{>0}$, coincides with 
the totally positive part of $\mathcal R_{v,w}$ defined 
by Lusztig \cite{Lus:TotPos94}, see  \cite[Theorem~11.3]{MarRie:ansatz}. 
\end{rem}

\section{The mirror conjecture for $G/P$}\label{s:conjecture}
 
Let $Z^{\mathfrak h^\vee}_P
\to (\mathfrak h^\vee)^{W_P}$ be the 
pullback of the family $pr_1:Z_P\to (T^\vee)^{W_P}$ 
under the exponential map $\exp:(\mathfrak h^\vee)^{W_P}\to (T^\vee)^{W_P}$. So, explicitly,
\begin{equation*}
Z^{\mathfrak h^\vee}_P=\{(h^\vee,b)\in (\mathfrak h^\vee)^{W_P}\x B^\vee_-\ | \
 b\in U^\vee_+\exp(h^\vee)\dot w_P\dot w_0\inv U^\vee_+\ \}.
\end{equation*}
For $h^\vee$ in $\mathfrak h^\vee$ we write $Z_P^{h^\vee}$ for the fiber  over $h^\vee$
in $Z^{\mathfrak h^\vee}_P$.
We may identify this fiber with
$$
B^\vee_-\cap U^\vee_+\exp(h^\vee)\dot w_P\dot w_0\inv
U^\vee_+.$$
Note that as in Section~\ref{s:mirrors}, 
\begin{eqnarray*}
Z^{\mathfrak h^\vee}_P&\overset\sim\To &(\mathfrak h^\vee)^{W_P}\x\mathcal R^\vee_{w_P,w_0},\\
(h^\vee,b) & \mapsto & (h^\vee,b\dot w_0 B^\vee_-).
\end{eqnarray*}
The phase function $\mathcal F_P$ 
pulled back to $Z^{\mathfrak h^\vee}_P$ 
will be again denoted by $\mathcal F_P$. 

Now let $\mathbf i_0$ be a reduced expression of $w_0$ and  $\omega$ 
the $n_P$-form on $\mathcal R^\vee_{w_P,w_0}$ defined in 
Proposition~\ref{p:omega}. Let us pull this
$n_P$-form back to $Z^{\mathfrak h^\vee}_P$ by the map 
\begin{eqnarray*}
Z^{\mathfrak h^\vee}_P &\To &\mathcal R^\vee_{w_P,w_0},\\
(h^\vee,b) &\mapsto & b\dot w_0 B^\vee_-/B^\vee_-,
\end{eqnarray*}
and denote the resulting form again by 
$\omega$. Note that $\omega$ 
depends on the reduced expression $\mathbf i_0$ only 
for its sign. We write $\omega_{h^\vee}$ for the restriction of
$\omega$ to the fiber $Z_P^{h^\vee}$.
 
\begin{conj}\label{c:mirrorconjecture}
The integrals \eqref{e:S} defined in terms of the 
mirror datum $(Z_P^{\mathfrak h^\vee},\omega,\mathcal F_P)$
give solutions to the  quantum 
differential equations \cite{Giv:EquivGW,CoxKatz:QCohBook} of $G/P$.
\end{conj}

We now want to state a $T$-equivariant version of the above conjecture.
For this we need to integrate over functions defined on the covering $\tilde Z_P$
of $Z_P$. We therefore pull back also this covering family 
$pr_1:\tilde Z_P\to (T^\vee)^{W_P}$ to $(\mathfrak h^\vee)^{W_P}$,
to get
$$
\tilde Z_P^{\mathfrak h^\vee}=\{(h^\vee,b,h_R)\ |\  (\exp(h^\vee),b,h_R)\in \tilde Z_P\}.
$$ 
The pullbacks of $\tilde\phi$ and $\mathcal F_P$ to $\tilde Z_P^{\mathfrak h^\vee}$
will again be denoted by $\tilde\phi$ and $\mathcal F_P$, respectively. Moreover
the map $\tilde Z_P^{\mathfrak h^\vee}\to Z_P^{\mathfrak h^\vee}$ which forgets $h_R$ is again a covering, and the
$n_P$-form $\omega$ on $Z_P^{\mathfrak h^\vee}$ pulls back to an $n_P$-form on 
$\tilde Z_P^{\mathfrak h^\vee}$ which we denote by $\tilde \omega$. 
The restriction of $\tilde\omega$ to a fiber  $\tilde Z_P^{h^\vee}$ 
of the family $pr_1:\tilde Z_P^{\mathfrak h^\vee}\to(\mathfrak h^\vee)^{W_P}$  is denoted by
$\tilde\omega_{h^\vee}$.
\begin{conj}\label{c:EqMirrorconjecture}
A full set of solutions to the $T$-equivariant quantum 
differential equations \cite{Giv:EquivGW,CoxKatz:QCohBook} of $G/P$ 
is given by integrals
\begin{equation}\label{e:Iagain}
\tilde S_{\tilde\Gamma}(h^\vee,h)=\int_{\tilde\Gamma_{h^\vee}}e^{\mathcal F_P/\hbar}\tilde\phi(\ ,h)\ \tilde \omega_{h^\vee},
\end{equation}
where $\tilde \Gamma=(\tilde\Gamma_{h^\vee})_{h^\vee\in\mathfrak h_{\R}^\vee} $ is a
continuous family of suitable integration contours $\tilde\Gamma_{h^\vee}$ in
$\tilde Z_P^{h^\vee}$ such that $\tilde S_{\tilde\Gamma}$ converges. 
\end{conj}

We note here that the equivariant quantum cohomology ring of $G/P$ is 
semisimple. (This follows
from the fact that equivariant cohomology is semisimple.) 
Correspondingly, in a generic fiber determined by an $h^\vee\in (\mathfrak h^\vee)^{W_P}$ and for 
generic $h\in \mathfrak h$, the function 
$\mathcal F_{P}+\ln \phi(\  ,h)$ has the correct number, $\dim H^*(G/P)$, of 
non-degenerate critical points counted in $Z_P^{h^\vee}$. This suggests that one
could be able to construct the right number of suitable integration contours
using Morse theory,
as asserted  in the $SL_{n+1}/B$ case by Givental, and Joe and Kim,
which would hopefully give a basis of solutions. The same need not be true for general 
$G/P$ in the  
non-equivariant setting, that is for $h=0$ above.  However, the non-equivariant quantum 
cohomology ring is also known to be always semisimple for the full flag variety $G/B$ \cite{Kos:QCoh}, 
and for Grassmannians, 
see \cite{Gepner:FusRing, SiTi:QCoh}. 
 
\section{The mirror constructions for $SL_{n+1}/B$ of Givental, Joe and Kim}
  
Givental  constructed a mirror family for $SL_{n+1}/B$ in \cite{Giv:MSFlag}. In this section we recall Givental's construction and identify his mirror family with a restriction of ours to an open subset. 
We will show in that case 
that the oscillatory integrals \eqref{e:S} arising from our mirror construction 
agree with those of Givental, proving Conjecture~\ref{c:mirrorconjecture} in that case.
In the $T$-equivariant setting an analogue of Givental's
mirror theorem was given by Joe and Kim \cite{JoeKim:EquivMirrors}. We go on to
review their construction and compare it with ours for the equivariant case, showing 
that  the integrals constructed by Joe and Kim can indeed be written in the form \eqref{e:Iagain}. 
This supports our Conjecture~\ref{c:EqMirrorconjecture}.

We note also that Batyrev, Ciocan-Fontanine, Kim and van Straten \cite{BCKS:MSPFlag}  
proposed a mirror family for $SL_{n+1}/P$ in the style of Givental's. For the direct 
relationship between their construction and the type $A$ Peterson variety
see \cite{Rie:TotPosGBCKS}. 

\subsection{} Let  $G=SL_{n+1}$, so $G^\vee=PSL_{n+1}$. We use the standard 
choice of Chevalley generators
$e_i=e_i^\vee=E_{i,i+1}$ and $f_i=f_i^\vee=E_{i+1,i}$, where $E_{j,k}$ is the matrix with $1$ in position $(j,k)$
and zeros elsewhere. Correspondingly we have the simple root subgroups $x_i(t)=\mathbf 1_{n+1}+
t E_{i,i+1}$ which we may consider to be lying in $SL_{n+1}$ or $PSL_{n+1}$, depending on the context.
We now recall the type $A$ mirror construction from \cite{Giv:MSFlag}.
\subsubsection{}\label{s:quiver}
Consider the
quiver $(\mathcal V,\mathcal A)$ which looks as follows
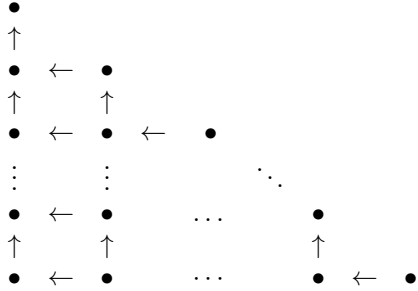
\begin{equation*}
\begin{matrix}
\bullet    &             &           &&           &&&&&\\
\uparrow &             &           &&           &&&&&\\
\bullet    & \leftarrow &\bullet    &&           &&&&&\\
\uparrow &             &\uparrow &&           &&&&&\\
\bullet    & \leftarrow &\bullet    &\leftarrow & \bullet &&&&&\\
\vdots     &             & \vdots    &            &         &\ddots &&&&\\
\bullet &\leftarrow &\bullet &&\dots&&\bullet&&           \\
\uparrow &             & \uparrow          &            & &
           &\uparrow&&&\\
\bullet    &\leftarrow  &  \bullet  &              & \cdots &
&\bullet &\leftarrow
           &\bullet
\end{matrix}
\end{equation*}
We divide the set of vertices $\mathcal V$ up  into the vertices along the 
diagonal, $\mathcal
V_\circ=\{v_{11},\dotsc, v_{n+1,n+1}\}$, and the vertices below the diagonal, 
$\mathcal
V_-=\{v_{ij}\ |\  1\le j<i\le n+1\}$. The labeling is as for 
the entries of a matrix. Let 
$\mathcal A=\mathcal A_c\sqcup \mathcal A_d$ be the set of arrows, 
divided into vertical and horizontal arrows, respectively. For any arrow $a$
denote by $h_a$ and $t_a\in\mathcal V$ the head and tail of $a$. In fact let us
label the arrow $a$ by $c_{ij}$ if $a$ is a vertical arrow  and $h_a=v_{ij}$, and 
by $d_{ij}$ if $a$ is horizontal with $t_a=v_{ij}$.  
Let
\begin{equation*}
\mathcal Z:=\left\{\sigma=(\sigma_a)_{a\in \mathcal A}\in (\C^*)^{\mathcal A}\ |\
\sigma_d \sigma_c=\sigma_{c'}\sigma_{d'}\ \text{for all configurations
\eqref{e:box}}\ \right \},
\end{equation*}
where
\begin{equation}\label{e:box}
\begin{CD}\bullet &@<d<<&\bullet \\
@A{c'}AA & &@A AcA\\
\bullet &@<{d'}<< & \bullet
\end{CD}
\end{equation}
is a square in the quiver.

For simplicity of notation we identify the arrows with coordinate functions
on $\mathcal Z$. In other words we may think of $\mathcal A\subset \C[\mathcal Z]$
as invertible generators for the coordinate ring of $\mathcal Z$. Define 
\begin{equation*}
\tilde q_i=c_{ii} d_{i+1,i+1}\quad \in \C[\mathcal Z],
\end{equation*}
for $i=1,\dotsc, n$. Then we have a family of varieties 
\begin{equation}\label{e:family}
\tilde q=(\tilde q_1,\dots,\tilde q_{n}):\mathcal Z\To (\C^*)^{n}.
\end{equation}
Let the fiber over $Q\in(\C^*)^{n}$ be denoted by $\mathcal Z_Q$. 
This map is a trivial fibration with fiber isomorphic to
$(\C^*)^{\left({n+1}\atop 2\right)}$. 
Explicitly, consider the isomorphism
\begin{equation}
(\C^*)^{\mathcal V_-}\times (\C^*)^{n}
\overset\sim\To \mathcal Z
\end{equation}
given by $((z_v)_{v\in\mathcal V_-}\times (z_{v_{ii}})^n_{i=1})
\mapsto \sigma:=(z_{h_a}z_{t_a}\inv)_{a\in\mathcal A}$, where we 
set $z_{v_{n+1,n+1}}=1$. 
In particular we have  vertex coordinates 
given by $t_{v}(\sigma)=z_v$.
We will denote $t_{v_{ij}}$ also by $t_{ij}$ for convenience.
Note that $t_{{i,i}}
t_{{i+1,i+1}}\inv=\tilde q_i$. 

\subsubsection{}
The phase function on this family is 
defined to be
\begin{equation}\label{e:GivF}
\tilde {\mathcal F}(\sigma)=\sum_{a\in \mathcal A }
\sigma_a.
\end{equation}
One has the 
following simple description of the critical points 
of $\tilde {\mathcal F}$ along the fibers of $\tilde q$,  
\begin{equation}\label{e:GivFcrit}
\mathcal Z^{crit}=
\left \{\, \sigma\in \mathcal Z\ \left |\
\sum_{a\in\mathcal A,\ h_a=v}\sigma_a\ = \sum_{a\in\mathcal A, \
t_a=v}\sigma_a, \ \text{ for all}\ v\in\mathcal V_-\right.\right \}.
\end{equation}
\subsubsection{}
We will now recall the construction of the quantum Toda lattice
in type $A$ and state Givental's mirror theorem. Consider the map
\begin{eqnarray*}
\varepsilon:\C^{n+1} &\to &(\C^*)^n\\
(T_{i})_{i=1}^{n+1}&\mapsto & (\exp(T_i-T_{i+1}))_{i=1}^n.
\end{eqnarray*}
Let $\tilde {\mathcal Z}\to \C^{n+1}$ be the pullback of
the bundle $\tilde q:\mathcal Z\to {(\C^*)}^n$ 
by $\varepsilon$. We denote $\varepsilon^*(\tilde q)$ and  
$\varepsilon^*(\tilde{\mathcal F})$ again by $\tilde q$ and 
$\tilde{\mathcal F}$, respectively.

Solving the $\mathfrak {gl}_{n+1}$ quantum Toda lattice 
means finding smooth functions $S=S(T_1,\dotsc, T_{n+1})$
satisfying
\begin{equation}\label{e:TodaA}\det\left(
x \mathbf 1_{n+1}+ 
\begin{pmatrix}
\hbar\frac{\partial}{\partial T_1} & e^{T_1-T_2} &   &  & \\
- 1  &\hbar\frac{\partial}{\partial T_2} &e^{T_2-T_3}&& \\
   & -1 & \ddots  &\ddots&\\
    &     &      \ddots  & &e^{T_{n}-T_{n+1}} \\
     &     &            & -1 &\hbar\frac{\partial}{\partial T_{n+1}}
\end{pmatrix}\right )S=x^{n+1} S,
\end{equation}
where $\hbar\in\R_{>0}$.   
Note that the coefficients of the polynomial in $x$ on the left hand side are well-defined
differential operators. 
By \cite{Giv:MSFlag, Kim:QCohG/B} the  quantum differential equations for
$SL_{n+1}/B$ make up the quantum Toda lattice
for $\mathfrak {sl}_{n+1}$, whose solutions 
are obtained by restricting the solutions $S$ of \eqref{e:TodaA} to the 
subspace of $\C^{n+1}$ defined by $\sum_{i=1}^{n+1} T_i=0$.

\begin{thm}[Givental \cite{Giv:MSFlag}] \label{t:Givental}
Let  $\Gamma=(\Gamma_{T_*})_{T_*\in\C^{n+1}}$ be 
a continuous family of possibly non-compact  $\left( n+1\atop 2\right)$-cycles 
$$ 
\Gamma_{T_{*}}\subset \mathcal Z_{\varepsilon(T_1,\dotsc, T_{n+1})}
$$ 
obtained from descending Morse cycles for $Re(\mathcal F)$.
The integrals 
\begin{equation*}
\tilde S_\Gamma(T_1,\dotsc,T_{n+1})=\int_{\Gamma_{T_*}} e^{{\tilde{\mathcal F}}/{\hbar}}
 \bigwedge_{v\in \mathcal V_-} 
\frac{dt_{v}}{t_v}
\end{equation*}
solve the system of differential equations \eqref{e:TodaA}.  
\end{thm}   

\subsection{}
Let us now consider $G^\vee=PSL_{n+1}(\C)$.  
We have $n_B=\frac{n(n+1)}2$ and let $\mathbf
i_0=(i_1,\dotsc,i_{n_B})$ be the reduced expression 
of $w_0$ obtained by successively concatenating
the sequences $l_k=(n,\dotsc, k)$ for $k=1,\dotsc, n$. 
To $\sigma\in\mathcal Z$ we associate two unipotent upper-triangular matrices,
\begin{align}
x_{\mathbf c}(\sigma)&=\prod_{k=1}^{n}\left(\prod_{j=n}^k x_{j}(\sigma_{c_{j,k}})\right),\\
x_{\mathbf d}(\sigma)&=\prod_{k=n}^{1}\left(\prod_{j=1}^{k}
x_{n-j+1}(\sigma_{d_{k+1,j+1}})\right).
\end{align}
We also associate to $\sigma$ an element $\tau(\sigma)$ of the maximal 
torus $T^\vee$ of $PSL_{n+1}$, which is given by 
 $$
\begin{bmatrix}t_{11}(\sigma)&&&&\\
& t_{22}(\sigma) &&&\\
&&   \ddots&&\\
&  && t_{nn} (\sigma)&\\
& &        && 1\\
\end{bmatrix}
=
\begin{bmatrix}\tilde q_1\cdots \tilde q_{n}(\sigma)&&&&\\
& \ddots &&&\\
&&   \tilde q_{n-1} \tilde q_n(\sigma)&&\\
&  && \tilde q_n (\sigma)&\\
& &        && 1\\
\end{bmatrix}.
$$
The definition of a matrix $x_{\mathbf c}(\sigma)$ associated to a point in Givental's 
$\mathcal Z$ 
can already be found in \cite{GKLO:GaussGivental}
and \cite{Rie:TotPosGBCKS}. Its combination with 
$x_{\mathbf d}(\sigma)$ and $\tau(\sigma)$ required
for the full Lie theoretic interpretation of $\mathcal Z$, see the theorem below, 
appears here for the first time.

\begin{thm}\label{t:compare} Let $Z_B$ and $\mathcal F_B$ 
be as defined in Section~\ref{s:mirrors}. 
\begin{enumerate} 
\item The map
\begin{eqnarray*}
\beta:\mathcal Z &\to & G^\vee \\
\sigma &\mapsto & x_{\mathbf c}(\sigma) \tau(\sigma)\dot w_0\inv x_{\mathbf d}(-\sigma)\inv
\end{eqnarray*}
has image in $B^\vee_-$ and, taken together with $\tau$, defines an open embedding 
\begin{equation*}
(\tau,\beta):\mathcal Z \hookrightarrow Z_B.
\end{equation*}
\item
We have
$(\tau,\beta)^*(\mathcal F_B)=\tilde{\mathcal F}$. 
In particular the map $\sigma\mapsto x_{\mathbf c}(\sigma) B^\vee_-/B^\vee_-$ identifies $\mathcal Z^{crit}$
with the intersection of $Y_B$ and the open Deodhar stratum 
in $\mathcal R^\vee_{1,w_0}$ corresponding to $\mathbf i_0$.
\item
For $Q\in (\C^*)^n$ let $\gamma_Q:\mathcal Z_Q\to \mathcal R^\vee_{1,w_0}$ be the 
embedding of a fiber given by 
$\gamma_Q(\sigma)= x_{\mathbf c}(\sigma) B^\vee_-/B^\vee_-$. Suppose $\omega$ is a holomorphic $n_B$-form
on $\mathcal R^\vee_{1,w_0}$ as in Proposition~\ref{p:omega}. Then
\begin{equation}
\gamma_Q^*(\omega)=\epsilon\bigwedge_{v\in\mathcal V_{-}}\frac{dt_v}{t_v},
\end{equation}
where $\epsilon\in\{\pm 1\}$ and is independent of $Q$.  
Here the $t_v$ are the vertex coordinates defined in \ref{s:quiver}.
\end{enumerate}
\end{thm}

Note that Theorem~\ref{t:compare} together with 
Givental's Theorem~\ref{t:Givental} implies the  
Conjecture~\ref{c:mirrorconjecture} for $SL_{n+1}/B$.

\begin{lem}\label{l:brr}
Let ${\bf t}=(\bar t_{v})_{v\in\mathcal V}\in (\C^*)^\mathcal V$, where we may write 
$\bar t_{ij}$ for $\bar t_{v_{ij}}$ for short. Let $b({\bf t})\in GL_{n+1}$ be defined
by  
$$
b({\bf t}):=x_\mathbf c((\bar  t_{h_a}\bar t_{t_a}\inv)_{a\in\mathcal A})
\begin{pmatrix}
\bar t_{11}&&&\\
& \bar t_{22} &&\\
&  &  \ddots &\\
&         &&\bar  t_{n+1,n+1}\\
\end{pmatrix}
\dot w_0\inv x_{\mathbf d}((-\bar t_{h_a}\bar t_{t_a}\inv)_{a\in\mathcal A})\inv.
$$
Then for the fundamental representation $V(\omega_k)$ we have 
\begin{equation*}
\left<b({\bf t})\cdot\, v^+_{\omega_k}\ ,\
v^+_{\omega_k}\right>= \left(\prod_{i=1}^{n+1}\bar t_{ii}\right)
 \left(\prod_{i-j=k}\bar t_{ij}\inv \right).
\end{equation*}
\end{lem}

\begin{proof}
Let $\{v_1,\dotsc, v_{n+1}\}$ be the standard basis of $\C^{n+1}$, and choose the standard
highest weight vector
$v^+_{\omega_k}:=v_1\wedge\dotsc\wedge v_k$
in $V(\omega_k)=\bigwedge^k \C^{n+1}$. Then we have the lowest weight vector 
$$
v_{\omega_k}^-:=\dot w_0\inv\cdot (v_1\wedge\dotsc\wedge v_k) = v_{n-k+2}\wedge\dotsc\wedge v_{n+1},
$$ 
and
\begin{multline}\label{e:tildebeval}
b({\bf t})\cdot v^+_{\omega_k}=x_{\mathbf c}\left( (\bar t_{h_a}\bar t_{ t_a}\inv)_{a\in\mathcal A}\right)
\begin{pmatrix} \bar t_{11}&&&\\
& \bar t_{22} &&\\
&&   \ddots&\\
&  && \bar t_{n+1,n+1}
\end{pmatrix}
\cdot
v_{n-k+2}\wedge\dotsc\wedge v_{n+1}\\
=
\left( \prod_{j=n-k+2}^{n+1}{\bar t_{jj}}\right) \ x_{\mathbf c}\left( (\bar t_{h_a}\bar t_{ t_a}\inv)_{a\in\mathcal A}\right)
\cdot
(v_{n-k+2}\wedge\dotsc\wedge v_{n+1}).
\end{multline}
Now note that, written out, 
\begin{align*}
x_{\mathbf c}((\sigma_a)_{a\in \mathcal A})=&x_n(\sigma_{c_{n,1}})x_{n-1}(\sigma_{c_{n-1,1}})\cdots\cdots x_{k+1}(\sigma_{c_{k+1,1}})\underline{x_{k}(\sigma_{c_{k,1}})\cdots\cdots x_2(\sigma_{c_{21}})x_{1}(\sigma_{c_{11}})}\\
&x_n(\sigma_{c_{n,2}})\cdots\cdots\cdots x_{k+2}(\sigma_{c_{k+2,2}}) \underline{x_{k+1}(\sigma_{c_{k+1,2}})x_{k}(\sigma_{c_{k,2}})\cdots\cdots x_{2}(\sigma_{c_{22}})}\\
&\vdots\\
&x_n(\sigma_{c_{n,n-k}})\underline{x_{n-1}(\sigma_{c_{n-1,n-k}})x_{n-2}(\sigma_{c_{n-2,n-k}})\cdots\cdots x_{n-k}(\sigma_{c_{n-k,n-k}})}\\
&\underline{x_n(\sigma_{c_{n,n-k+1}})x_{n-1}(\sigma_{c_{n-1,n-k+1}})\cdots\cdots x_{n-k+1}(\sigma_{c_{n-k+1,n-k+1}})}\\
&x_n(\sigma_{c_{n,n-k+2}})x_{n-1}(\sigma_{c_{n-1,n-k+2}})\cdots x_{n-k+2}(\sigma_{c_{n-k+2,n-k+2}})\\
&\vdots\\
&x_n(\sigma_{c_{n,n-1}})x_{n-1}(\sigma_{c_{n-1,n-1}})\\
&x_n(\sigma_{c_{n,n}}).
\end{align*}
Each $x_j(a)=\exp(a e_j)$ simply acts by $1+a e_j$ on $\bigwedge^k \C^{n+1}$, and it is not
hard to check that in order to get from the lowest to 
the highest weight space via $x_{\mathbf c}((\sigma_a)_{a\in \mathcal A})$ we need to take the $e_j$-summand precisely
from each of the underlined $x_j(\sigma_{c_a})$ factors. Since in our case
$\sigma_a=\bar t_{h_a}\bar t_{ t_a}\inv$, we have
$$
\sigma_{c_{j+k-1,j}}\sigma_{c_{j+k-2,j}}\dotsc \sigma_{c_{j,j}}=\bar t_{jj}\bar t_{j+k,j}\inv,
$$
for the resulting contribution of the $j$-th row above, and so we find that in total
$$
\left< x_{\mathbf c}\left( (\bar t_{h_a}\bar t_{t_a}\inv)_{a\in\mathcal A}\right)
\cdot
v_{\omega_k}^-,v_{\omega_k}^+\right>=
\prod_{j=1}^{n-k+1} \bar t_{jj}\bar t_{j+k,j}\inv=\left(\prod_{j=1}^{n-k+1}\bar t_{jj}\right)\left(\prod_{i-j=k}\bar t_{ij}\inv\right).  
$$
Combining this with \eqref{e:tildebeval} we 
see that 
$$
\left<\tilde b\cdot v_{\omega_k}^+,v_{\omega_k}^+\right>=\left(\prod_{j=1}^{n+1}\bar t_{jj}\right)\left(\prod_{i-j=k}\bar t_{ij}\inv\right)
$$ 
and the lemma is proved.
\end{proof}
\vskip .2cm
\noindent{\it Proof of Theorem~\ref{t:compare}.~}
We assume for the moment that we have proved that $\beta$ has image in $B^\vee_-$,
and observe how the rest of the theorem follows from this assertion.
If $\beta$ has image in $B^\vee_-$, then $(\tau,\beta)$ defines a map  $\mathcal Z\to Z_B$. 
From this point of view the $c_{ij}$ 
correspond 
precisely the standard coordinates for $x_{\mathbf c}(\sigma)B^\vee_-/B^\vee_-=\beta(\sigma)\dot w_0
B^\vee_-/B^\vee_-$
in the open Deodhar stratum $\mathcal R^\vee_{\mathbf 1_+(\mathbf i_0),\mathbf i_0}$. 
Moreover the values of the $c_{ij}$ together with those of the $\tilde q_i$ suffice to determine
a point in $\mathcal Z$ uniquely.  Therefore we see that $(\tau,\beta)$ is injective, and its image
is equal
to the preimage of $T^\vee\x\mathcal R^\vee_{\mathbf 1_+(\mathbf i_0),\mathbf i_0}$ under the
trivialization \eqref{e:triv} of $Z_B$. This proves (1). 
Part (2) then follows from (1) and Lemma~\ref{l:formula} combined with Theorem~\ref{t:main}. See also 
\cite{Rie:TotPosGBCKS}. The third part of the theorem is an easy consequence of the definition of 
$x_{\mathbf c}$.

It remains to prove that $\beta$ has image in $B^\vee_-$.  We may
work in $GL_{n+1}$, rather than $PSL_{n+1}$, choosing the 
representative for $\tau(\sigma)$ as the one from its definition. 
Multiplying out the product $x_{\mathbf c}$ we obtain a matrix
\begin{equation*}
x_{\mathbf c}=\left(\begin{array}{ccccc}
 1 & G^{(1)}_1 & G^{(2)}_2& \dots  & G^{(n)}_{n}\\
  &  1        & G^{(2)}_1&        & G^{(n)}_{n-1}\\
   &           & 1      &        & \vdots        \\
   &           &           & \ddots    &  G^{(n)}_1 \\
  &           &           &          &       1
 \end{array}\right)
\end{equation*}
with entries given by
$$
G^{(j)}_k=\sum_{1\le m_1<\cdots<m_k\le j}\left(
\prod_{i=1}^k {c_{j-k+i,m_i}}\right).
$$ 
Similarly let $\tilde x(\sigma):=x_{\mathbf d}(-\sigma)\inv$. Then
\begin{equation*}
\tilde x=\left(\begin{array}{ccccc}
 1 & \tilde G^{(1)}_1& \tilde G^{(2)}_2& \dots  & \tilde G^{(n)}_{n}\\
  &  1        & \tilde G^{(2)}_1&        & \tilde G^{(n)}_{n-1}\\
   &           & 1      &        & \vdots        \\
   &           &           & \ddots    & \tilde  G^{(n)}_1\\
  &           &           &          &       1
 \end{array}\right)
\end{equation*}
where 
\begin{equation*}
\tilde G^{(j)}_k=\sum_{n-j+2\le m_1\le\dotsc\le m_k\le n+1}\left( \prod_{i=1}^k
d_{m_{k-i+1},n+1-j+i}\right).
\end{equation*}

The $(j,r+1)$ entry of the matrix $\beta(\sigma)=x_{\mathbf c}(\sigma)\tau(\sigma)\dot w_0\inv x_{\mathbf d}(-\sigma)\inv$ is
\begin{equation}\label{e:entry} \beta_{j,r+1}:=
(0,\dotsc,0,1,G^{(j)}_1,G^{(j+1)}_2,\dotsc, G^{(n)}_{n-j+1})\cdot
\begin{pmatrix}
0\\
\vdots \\
0\\
\pm \tilde q_{n-r+1}\cdots \tilde q_n 1\\
\vdots\\
\tilde q_{n-1}\tilde  q_{n} \tilde G^{(r)}_{r-2}\\
-\tilde q_n\tilde G^{(r)}_{r-1}\\
\tilde G^{(r)}_r
\end{pmatrix}
\end{equation} 
evaluated at $\sigma$. We want to show that this expression is zero
when $j\le r$. 

In the rank $1$ case we have for $j=r=1$ 
\begin{equation*}
(1,G^{(1)}_1)\cdot 
\begin{pmatrix}
-\tilde q_1\\
\tilde G^{(1)}_1
\end{pmatrix}=-\tilde q_1 + c_{11}d_{22}=0.
\end{equation*}
We will prove the general case by induction. 

Consider the two embeddings of $GL_n$ into $GL_{n+1}$ 
corresponding to the  subsets $I_L=\{1,\dotsc, n-1\}$ and
$I_R=\{2,\dotsc, n\}$ of $I$. The first gives the subgraph 
$(\mathcal V_L,\mathcal A_L)$ of $(\mathcal V,\mathcal A)$
obtained by erasing the last row of vertices. And the second
gives the subgraph $(\mathcal V_R,\mathcal A_R)$ where the 
first column has been removed.

We add superscripts $L$ and $R$ to any of the matrices
 $\tilde x,  x_{\mathbf c},\tau,\beta$ if we 
 are referring to their analogues defined in terms of the 
 graphs $(\mathcal V_L,\mathcal A_L)$ or $(\mathcal V_R,\mathcal A_R)$,
 respectively.  We denote by 
 $\tilde G^{(r,L)}_{k}$ and $\tilde G^{(r,R)}_k$ the matrix coefficients
 of $\tilde x^L$ and $\tilde x^R$, each viewed inside its respective 
 copy of $GL_n$. Similarly for $x_{\mathbf c}^L$ and $x_{\mathbf c}^R$
 and their entries. 

It is easy to check that   
\begin{equation*}
\tilde G^{(r)}_k=\tilde G^{(r-1,L)}_k+d_{n+1,n-r+2}\tilde G^{(r-1,R)}_{k-1}.
\end{equation*}
So we have
\begin{equation}\label{e:summands}
\begin{pmatrix}
0\\
\vdots \\
0\\
\pm \tilde q_{n-r+1}\dotsc \tilde q_n 1\\
\mp \tilde q_{n-r+2}\dotsc\tilde q_n \tilde G^{(r)}_{1}\\
\vdots\\
\tilde q_{n-1}\tilde  q_{n} \tilde G^{(r)}_{r-2}\\
-\tilde q_n\tilde G^{(r)}_{r-1}\\
\tilde G^{(r)}_r
\end{pmatrix}
=
\begin{pmatrix}
0\\
\vdots \\
0\\
\pm \tilde q_{n-r+1}\dotsc \tilde q_n 1\\
\mp \tilde q_{n-r+2}\dotsc\tilde q_n \tilde G^{(r-1,L)}_{1}\\
\vdots\\
\tilde q_{n-1}\tilde  q_{n} \tilde G^{(r-1,L)}_{r-2}\\
-\tilde q_n\tilde G^{(r-1,L)}_{r-1}\\
0
\end{pmatrix}
+
d_{n+1,n-r+2}
\begin{pmatrix}
0\\
\vdots\\
\vdots \\
0\\
\mp \tilde q_{n-r+2}\dotsc \tilde q_n 1\\
\vdots\\
\tilde q_{n-1}\tilde  q_{n} \tilde G^{(r-1,R)}_{r-3}\\
-\tilde q_n\tilde G^{(r-1,R)}_{r-2}\\
\tilde G^{(r-1,R)}_{r-1}
\end{pmatrix}.
\end{equation} 
We now want to evaluate \eqref{e:entry} using \eqref{e:summands}. The first summand gives
a contribution of 
\begin{equation*}
\left(0,\dotsc,0,1,G^{(j)}_1,G^{(j+1)}_2,\dotsc, G^{(n)}_{n-j+1}\right )\cdot
\begin{pmatrix}
0\\
\vdots \\
0\\
\pm \tilde q_{n-r+1}\dotsc \tilde q_n 1\\
\mp \tilde q_{n-r+2}\dotsc\tilde q_n \tilde G^{(r-1,L)}_{1}\\
\vdots\\
\tilde q_{n-1}\tilde  q_{n} \tilde G^{(r-1,L)}_{r-2}\\
-\tilde q_n\tilde G^{(r-1,L)}_{r-1}\\
0
\end{pmatrix}
\end{equation*}
to \eqref{e:entry}. This equals $-\beta_{j,r}^L$ and is therefore zero if $j<r$ by the induction 
hypothesis, where we are 
reducing to the graph with the last row removed. Note that in this induction step 
we are also removing the
last row of $\tau$, which we had normalized to $1$ in the rank $n$ case. 
This accounts for the apparent factor of $\tilde  q_{n}$ in the above formula for $-\beta_{j,r}^L$.

For the second summand  notice that  we can decompose the 
entries of the row vector in \eqref{e:entry} using
\begin{equation*}
G^{(j+k)}_{k+1}=c_{j,1}G^{(j+k,R)}_{k}+ G^{(j+k,R)}_{k+1}.
\end{equation*}
Then we get $d_{n+1,n-r+2}(c_{j,1}\beta^R_{j,r} +\beta^R_{j,r+1})$, which is also zero whenever $j<r$,
 by the induction hypothesis
this time applied to the graph with the left most column removed. 

Thus we have seen that \eqref{e:entry} vanishes whenever $j<r$. 
If $j=r$ we are left with two nonzero summands, giving
\begin{equation}\label{e:rr}
\beta_{r,r+1}=-\beta^L_{r,r} + d_{n+1,n-r+2}c_{r,1}\beta^R_{r,r}.
\end{equation}  
It remains to show that this matrix coefficient vanishes.

By induction assumption $\beta^L$ and $\beta^R$ are upper-triangular,
so we have
\begin{equation}\label{e:bX}
\left<\beta^X\cdot\, v^+_{\omega_r}\ ,\
v^+_{\omega_r}\right>= \beta^X_{1,1}\beta^X_{2,2}\cdots \beta^X_{r,r},
\end{equation}
for $X=L$ or $R$. Now let $\sigma\in\mathcal Z$ and consider the vertex
coordinates $\bar t_{v}=t_{v}(\sigma)$. Then the corresponding 
`truncated' elements are $\sigma_L=(\bar t_{h_a}\bar t_{t_a}\inv)_{a\in\mathcal A_L}$
and $\sigma_R=(\bar t_{h_a}\bar t_{t_a}\inv)_{a\in\mathcal A_R}$,
and Lemma~\ref{l:brr} says that 
\begin{align*}
&\left<\beta^L(\sigma_L)\cdot\, v^+_{\omega_k}\ ,\
v^+_{\omega_k}\right>= \left(\prod_{i=1}^{n}\bar t_{ii}\right)
 \left(\prod_{j=1}^{n-k}{\bar t_{j+k,j}}\inv \right),\\
& \left<\beta^R(\sigma_R)\cdot\, v^+_{\omega_k}\ ,\
v^+_{\omega_k}\right>= \left(\prod_{i=2}^{n+1}\bar t_{ii}\right)
 \left(\prod_{j=2}^{n-k+1}\bar t_{j+k,j}\inv \right).
\end{align*}
Combining these formulas with \eqref{e:bX} we find that
\begin{align*}
&\beta^L_{r,r}(\sigma_L)= \left(\prod_{j=1}^{n-r+1}{\bar t_{j+r-1,j}} \right)\left(\prod_{j=1}^{n-r}{\bar t_{j+r,j}}\inv \right),
\\
&\beta^R_{r,r}(\sigma_R)= \left(\prod_{j=2}^{n-r+2}\bar t_{j+r-1,j} \right) \left(\prod_{j=2}^{n-r+1}\bar t_{j+r,j}\inv \right).
\end{align*}
Now substituting also $c_{r,1}(\sigma)=\bar t_{r,1}\bar t_{r+1,1}\inv$  and $ d_{n+1,n-r+2}(\sigma) =\bar t_{n+1,n-r+1}\bar t_{n+1,n-r+2}\inv$ it follows directly that
\begin{equation*}
d_{n+1,n-r+2}c_{r,1}\beta^R_{r,r}=\beta^L_{r,r}. 
\end{equation*}
This shows that $\beta_{r,r+1}=0$, by \eqref{e:rr}, and finishes the proof. 
\qed
\vskip .2cm

\subsection{The $T$-equivariant case}
In Joe and Kim's work \cite{JoeKim:EquivMirrors}, mirror symmetric
solutions to the $T$-equivariant quantum differential equations
of $SL_{n+1}/B$ are given as integrals over a function defined 
on a universal cover of $\mathcal Z$. We briefly review this construction
here and compare it with our definitions applied to the equivariant 
$SL_{n+1}/B$ case. 

The $T$-equivariant quantum differential equations are deformations
of the usual quantum differential equations by
the ring
$$
H^*_T( pt )=\C[\mathfrak h]=\C[\lambda_1,\dotsc, \lambda_{n+1}]/(\sum {\lambda_i}).
$$ 
Namely in the $T$-equivariant case the differential
equations to solve are obtained by replacing \eqref{e:TodaA} by 
\begin{equation}\label{e:EquivTodaA}\det\left(
x \mathbf 1_{n+1}+ 
\begin{pmatrix}
\hbar\frac{\partial}{\partial T_1} & e^{T_1-T_2} &   &  & \\
- 1  &\hbar\frac{\partial}{\partial T_2} &e^{T_2-T_3}&& \\
   & -1 & \ddots  &\ddots&\\
    &     &      \ddots  & &e^{T_{n}-T_{n+1}} \\
     &     &            & -1 &\hbar\frac{\partial}{\partial T_{n+1}}
\end{pmatrix}\right )\tilde S=\prod_{i=1}^{n+1} (x+\lambda_i) \tilde S.
\end{equation}

To generalize Givental's mirror theorem to solve \eqref{e:EquivTodaA}
Joe and Kim deform the phase function $\tilde{\mathcal F}$ of Givental, or
more precisely they first pull it back to a universal cover of $\mathcal Z$
and then deform it there.

\begin{defn}\label{d:Ztilde}
Recall the notations from Section~\ref{s:quiver}. We let 
$$
\tilde{\mathcal Z}:=\{(T_v)_{v\in \mathcal V}\in \C^{\mathcal V}\ |\ \sum T_{v_{ii}}=0\}.
$$
For given $(T_1,\dotsc, T_{n+1})\in\C^{n+1}$ with $\sum_i {T_i}=0$ we define
$$
\tilde{\mathcal  Z}_{(T_1,\dotsc, T_{n+1})}:=\{(T_v)_{v\in \mathcal V}\ |\ T_{v_{ii}}=T_i, \text{ for $i=1,\dotsc, n+1$}\}.
$$
\end{defn}

The map $c:\tilde{\mathcal  Z}\to \mathcal Z$ given by
$$
c:(T_v)_{v\in \mathcal V}\mapsto (e^{T_{h_a}-T_{t_a}})_{a\in\mathcal A}
$$
makes $\tilde{\mathcal  Z}$ into a universal covering space for $\mathcal Z$.
We may think of the $T_v$ as logarithmic vertex variables, although 
the $\exp(T_v)$ recover the vertex variables $t_v$ from 
Section~\ref{s:quiver} only up to a common scalar multiple, as we are 
working with a different normalization now~: We have 
$\sum T_{v_{ii}}=0$ rather than  $T_{v_{n+1,n+1}}=0$.

To deform Givental's phase function Joe and Kim attach 
`weights'  depending on the parameters $\lambda_i$
to the edges of the graph $(\mathcal V,\mathcal A)$ as
follows,
\begin{eqnarray*}
\lambda_{c_{i1}}&:= &\lambda_i+\frac 1 2(\lambda_1+\dotsc +\lambda_{i-1}),\\
\lambda_{c_{ij}}&:= &\frac 1 2\lambda_{i-j+1}, \qquad \text{if $j>1$,}\\
\lambda_{d_{n+1,j}}&:= &-\lambda_{n+1-j}-\frac 1 2(\lambda_1+\dotsc+\lambda_{n-j}),\\
\lambda_{d_{ij}}&:= &-\frac 1 2\lambda_{i-j+1}, \ \qquad\text{if $i<n+1$,}
\end{eqnarray*}
and set
\begin{equation}\label{e:FJK}
\tilde{\mathcal F}_{JK}\left((T_{v})_{v\in\mathcal V};(\lambda_i)_i\right):=
\tilde{\mathcal F}\left((e^{T_{h_a}-T_{t_a}})_{a\in \mathcal A}\right)+
\sum_{a\in\mathcal A} \lambda_a (T_{h_a}-T_{t_a}).
\end{equation}

\begin{thm}[Joe and Kim \cite{JoeKim:EquivMirrors}] \label{t:JoeKim}
Let $T_*$ run through the
$(T_1,\dotsc, T_{n+1})\in \C^{n+1}$ with $\sum T_i=0$, and let $\Gamma=(\Gamma_{T_*})_{T_*}$
be 
a continuous family of possibly non-compact  $\left( n+1\atop 2\right)$-cycles,
$$ 
\Gamma_{T_{*}}\subset \tilde{\mathcal  Z}_{(T_1,\dotsc, T_{n+1})},
$$ 
obtained as descending Morse cycles for $Re(\tilde{\mathcal F}_{JK})$.
The integrals 
\begin{equation*}
\tilde S_\Gamma(T_1,\dotsc,T_{n+1};\lambda_1,\dotsc,\lambda_{n+1})=\int_{\Gamma_{T_*}} e^{\frac 1\hbar {\tilde {\mathcal F}_{JK}}}
 \bigwedge_{v\in \mathcal V_-} 
dT_{v}
\end{equation*}
solve the $T$-equivariant quantum differential equations 
associated to $SL_{n+1}/B$.  
\end{thm}

We want to now give an explicit lift of the comparison map $(\tau,\beta):\mathcal Z\to Z_B$ 
and extend our Theorem~\ref{t:compare} about comparing phase functions to the equivariant case. 

\begin{defn}\label{d:gammaR}
Let $\gamma_R:\tilde{\mathcal  Z}\to \mathfrak h^\vee$ be defined by
$$
\omega_k^\vee(\gamma_R((T_v)_{v\in\mathcal V}))=-\sum_{ i-j=k}T_{v_{ij}}.
$$
Written out explicitly, $\gamma_R((T_v)_{v\in\mathcal V})$ is
$$
\begin{pmatrix}
-\underset{\scriptscriptstyle{i-j=1}}\sum T_{v_{ij}} & &&&\\
&(\underset{\scriptscriptstyle{i-j=1}}\sum T_{v_{ij}})-(\underset{\scriptscriptstyle{i-j=2}}\sum T_{v_{ij}}) &&&\\
&&\ddots &&\\
&&& (T_{v_{n1}}+T_{v_{n+1,2}})- T_{v_{n+1,1}}\\
&&&&T_{v_{n+1,1}}
\end{pmatrix}.
$$
We also define
\begin{equation*}
\tilde \beta:=\beta\circ c:\tilde{\mathcal  Z}\to B^\vee_-,
\end{equation*}
and a map
$
\gamma:\tilde{\mathcal  Z}\to \mathfrak h^\vee
$,
$$
\gamma((T_v)_{v\in\mathcal V})=
\begin{pmatrix}
T_{v_{11}}&&&\\
&T_{v_{22}}&&\\
&& \ddots &\\
&&& T_{v_{n+1,n+1}}
\end{pmatrix}.
$$
\end{defn}

\begin{thm}\label{t:EquivCompare}
\begin{enumerate}
\item
The maps $\gamma, \tilde \beta$ and $\gamma_R$ define a map
$$
(\gamma,\tilde\beta,\gamma_R):\tilde{\mathcal  Z}\to \tilde Z_B^{\mathfrak h^\vee},
$$
which is a covering composed with an open embedding.
Moreover $(\gamma,\tilde\beta,\gamma_R)$ takes the fiber $\tilde{\mathcal  Z}_{(T_1,\dotsc,T_{n+1})}$
to the fiber $\tilde Z^{h^\vee}_B$, where 
$h^\vee$ is the diagonal matrix 
with entries
$(T_1,\dotsc, T_{n+1})$.  
\item 
We have 
$$
(\gamma,\tilde\beta,\gamma_R)^*(\tilde\omega)=\pm\bigwedge_{v\in\mathcal V_-} dT_v
$$
for the pullback of our form $\tilde \omega$ from Section~\ref{s:conjecture} to 
$\tilde{\mathcal Z}$. 
\item
The integrand $e^{\tilde{\mathcal  F}_{JK}}$ of Joe and Kim
is obtained by pullback from our integrand,
$$
e^{\tilde{\mathcal  F}_{JK}}=(\gamma,\tilde \beta,\gamma_R;h)^*(e^{\mathcal F_B}\tilde \phi),
$$
where $h(\lambda_1,\dotsc,\lambda_{n+1})$ is the diagonal matrix with entries 
$\lambda_1,\dotsc, \lambda_{n+1}$. 
\end{enumerate}
\end{thm}

\begin{proof}
(1) For $(T_v)_{v\in\mathcal V}\in\tilde{\mathcal Z}$ we have that $\tilde{\beta}\left((T_v)_{v\in\mathcal V}\right)$ 
is given explicitly by
\begin{equation}\label{e:betatilde}
x_{\mathbf c}\left( (e^{T_{h_a}-T_{t_a}})_{a\in\mathcal A}\right)
\begin{bmatrix} e^{T_{v_{11}}}&&&\\
& e^{T_{v_{22}}} &&\\
&&   \ddots&\\
&  && e^{T_{v_{n+1,n+1}}}
\end{bmatrix}
\dot w_0\inv
x_{\mathbf d}\left( (-e^{T_{h_a}-T_{t_a}})_{a\in\mathcal A}\right)\inv.
\end{equation}
Here we have substituted $\sigma_a=e^{T_{h_a}-T_{t_a}}$ in the formula
from Theorem~\ref{t:compare}. Also note that the entries of the diagonal
matrix $\tau(\sigma)$ are $t_{ii}(\sigma)=e^{T_{v_{ii}}}/ e^{T_{v_{n+1,n+1}}}$ and,
as we are working in $PSL_{n+1}$, we can clear the denominators.  
From Theorem~\ref{t:compare} together with \eqref{e:betatilde} it is now immediate that 
$(\gamma,\tilde \beta)((T_v)_{v\in\mathcal V})\in Z^{\mathfrak h}_B$. 

To show that $(\gamma,\tilde \beta,\gamma_R)((T_v)_{v\in\mathcal V})$ lies in 
the covering space $\tilde Z^{\mathfrak h}_B$ it remains to prove that the diagonal part of $\tilde\beta\left((T_v)_{v\in\mathcal V}\right)$ is equal to $\exp(\gamma_R((T_v)_{v\in\mathcal V}))$.
For this let us consider the lower-triangular matrix $\tilde b$ in $SL_{n+1}$ which covers $\tilde\beta\left((T_v)_{v\in\mathcal V}\right)$
and is given by
\begin{multline}
x_{\mathbf c}\left( (e^{T_{h_a}-T_{t_a}})_{a\in\mathcal A}\right)
\begin{pmatrix} e^{T_{v_{11}}}&&&\\
& e^{T_{v_{22}}} &&\\
&&   \ddots&\\
&  && e^{T_{v_{n+1,n+1}}}
\end{pmatrix}
\dot w_0\inv
x_{\mathbf d}\left( (-e^{T_{h_a}-T_{t_a}})_{a\in\mathcal A}\right)\inv.
\end{multline}
Then by Lemma~\ref{l:brr} we have 
\begin{equation}
\left<\tilde b\cdot v^+_{\omega_k},v^+_{\omega_k}\right>=
e^{-\sum_{i-j=k}T_{v_{ij}}},
\end{equation}
using also that $\sum T_{v_{ii}}=0$. This implies the rest of (1),
comparing also with Definition~\ref{d:gammaR}.

\vskip .2cm

Part (2) of the theorem is  a consequence of Theorem~\ref{t:compare} (3),
using that $dT_v$ is the pullback 
to $\tilde {\mathcal Z}$ of $dt_v/t_v$.

\vskip .2cm

It now remains to show (3), namely that
$$
(e^{\mathcal F_B}\tilde\phi) \circ(\gamma,\tilde\beta, \gamma_R;h )
=e^{\tilde{\mathcal F}_{JK}}.
$$
By Theorem~\ref{t:compare} we already know 
that $e^{\mathcal F_B}\circ(\gamma,\tilde\beta, \gamma_R )=e^{\tilde{\mathcal F}\circ c}$. Therefore 
we only need to compare the effect of Joe and Kim's correction term with our factor $\tilde\phi$.
By definition
\begin{equation*}
\tilde\phi \circ(\gamma,\tilde\beta, \gamma_R;h )\ ((T_v)_{v\in\mathcal V};(\lambda_i)_i)=e^{<h((\lambda_i)_i)
,\gamma_R((T_v)_{v\in\mathcal V})>},
\end{equation*}
and the exponent evaluates to
\begin{equation}\label{e:exponent}
<h((\lambda_i)_i)
,\gamma_R((T_v)_{v\in\mathcal V})>=\sum_{k=1}^n(\lambda_{k+1}-\lambda_{k})\left(\sum_{i-j=k}T_{ij}\right).
\end{equation}
However, the weight factors of Joe and Kim are chosen precisely so that for every vertex $v=v_{ij}$
with $i-j=k$,
$$
\sum_{a, h_a=v} \lambda_a-\sum_{a, t_a=v} \lambda_a\ =\ \lambda_{k+1}-\lambda_{k}.
$$
Therefore Joe and Kim's correction term $\sum_{a\in\mathcal A} \lambda_a(T_{h_a}-T_{t_a})$,
reordered as a sum of $T_v$'s with $\lambda_i$ coefficients, gives precisely \eqref{e:exponent},
and we are done.  
\end{proof}
 
\subsection{}\label{s:MorseCycles}
To show that the solutions to the equivariant quantum differential equations
constructed by Joe and Kim can be put in the form of our 
Conjecture~\ref{c:EqMirrorconjecture}, we need finally to 
argue that our comparison map $(\gamma,\tilde\beta,\gamma_R)$
is one-to-one when restricted to the integration
contours put forward by Joe and Kim. 

Recall that $(\gamma,\tilde\beta, \gamma_R ):\tilde{\mathcal  Z}
\to \tilde Z_B^{\mathfrak h^\vee}$ defines a covering onto its image. 
Let us choose compatible Riemann metrics on 
$\tilde{\mathcal  Z}$ and $(\gamma,\tilde\beta, \gamma_R )(\tilde{\mathcal  Z})$, so that one is 
pulled back from the other.   
Suppose $p_0\in\tilde {\mathcal Z}_{T_{*}}$ is a critical point of $\tilde{\mathcal F}_{JK}$
with its corresponding `descending Morse cycle' for $Re(\tilde{\mathcal F}_{JK})$, denoted $\Gamma_{T_*}$. The gradient flow of $Re(\tilde{\mathcal F}_{JK})$ starting at 
$p\in\Gamma_{T_*}$ should therefore approach $p_0$ in the positive limit. 
This gradient flow maps out a curve which can also be obtained as the unique lifting 
through $p$ of the gradient 
flow curve of $Re(\mathcal F_B\tilde \phi)$ starting at $\bar p:=(\gamma,\tilde\beta, \gamma_R )(p)$
in the base.
Suppose now there was another point $p'$ in $\Gamma_{T^*}$ with 
the same image $\bar p':=(\gamma,\tilde\beta, \gamma_R )(p')=\bar p$. 
Then the gradient flow curve below, connecting $\bar p=\bar p'$ with the image $\bar p_0$ of the critical $p_0$, would have a lift through $p$ which ends up at $p_0$, and another lift through $p'$ which
also ends at $p_0$. This, however, is in contradiction with the unique lifting of curves property
of our covering. 

So we have seen that no two points in $\Gamma_{T^*}$ can map
to the same point 
under $(\gamma,\tilde\beta, \gamma_R) $. 
Therefore the map $(\gamma,\tilde\beta, \gamma_R)$ amounts to a change of coordinates
on the integration contour, and moreover by Theorem~\ref{t:EquivCompare}, a change of coordinates under which Joe and Kim's integrals 
transform to ones of the form \eqref{e:Iagain}.


\def\cprime{$'$}
\providecommand{\bysame}{\leavevmode\hbox to3em{\hrulefill}\thinspace}
\providecommand{\MR}{\relax\ifhmode\unskip\space\fi MR }
\providecommand{\MRhref}[2]{%
  \href{http://www.ams.org/mathscinet-getitem?mr=#1}{#2}
}
\providecommand{\href}[2]{#2}

\end{document}